\documentclass[11pt,dvipdfmx
]{amsart}
\usepackage{amsfonts}
\usepackage{amsmath, amsthm, amssymb,color}
\usepackage{latexsym}
\usepackage{txfonts}
\usepackage{bm}
\usepackage{graphicx}
\usepackage{graphics}
\usepackage{epsfig}

\numberwithin{equation}{section}
\allowdisplaybreaks

\newtheorem{lemma}{Lemma}[section]

\newtheorem{theorem}[lemma]{Theorem}
\newtheorem{proposition}[lemma]{Proposition}
\newtheorem{definition}[lemma]{Definition}
\newtheorem{corollary}[lemma]{Corollary}
\newtheorem{example}[lemma]{Example}
\newtheorem{exercise}[lemma]{Exercise}
\newtheorem{remark}[lemma]{Remark}
\newtheorem{fig}[lemma]{Figure}
\newtheorem{tab}[lemma]{Table}

\newcommand{\bth}{\begin{theorem}}
\newcommand{\ethe}{\end{theorem}}
\newcommand{\bre}{\begin{remark}\em }
\newcommand{\ere}{\end{remark}}
\newcommand{\ble}{\begin{lemma}}
\newcommand{\ele}{\end{lemma}}
\newcommand{\bde}{\begin{definition}}
\newcommand{\ede}{\end{definition}}
\newcommand{\bco}{\begin{corollary}}
\newcommand{\eco}{\end{corollary}}
\newcommand{\bpr}{\begin{proposition}}
\newcommand{\epr}{\end{proposition}}

\newcommand{\bexer}{\begin{exercise}}
\newcommand{\eexer}{\end{exercise}}

\newcommand{\bexam}{\begin{example}\rm  }
\newcommand{\eexam}{ \end{example}}

\newcommand{\bfi}{\begin{fig}}
\newcommand{\efi}{\end{fig}}

\newcommand{\btab}{\begin{tab}}
\newcommand{\etab}{\end{tab}}

\def\E{{\mathbb{E}}}
\def\P{{\mathbb{P}}}

\def\B_e{B_{\eta}(e)}

\renewcommand{\a}{\alpha }
\renewcommand{\b}{\beta}

\newcommand{\1}{{\bf 1}}

\newcommand{\wt}{\widetilde}

\newcommand{\blue}{\color{blue} }
\definecolor{darkblue}{rgb}{0,0,1}
\definecolor{darkgreen}{rgb}{0,1,0}
\definecolor{darkred}{rgb}{1, 0,0}

\newcommand{\garch}{{\rm GARCH}$(1,1)$}

\newcommand{\sta}{St\u aric\u a}

\newcommand{\asy}{asymptotic}

\newcommand{\ts}{time series}

\newcommand{\garchpq}{{\rm GARCH}$(p,q)$}

\newcommand{\bfu}{{\bf u}}

\newcommand{\bfR}{{\bf R}}
\newcommand{\bfC}{{\bf C}}

\newcommand{\sre}{stochastic recurrence equation}

\newcommand{\fidi}{finite-dimensional distribution}
\newcommand{\rv}{random variable}

\newcommand{\cov}{{\rm cov}}
\newcommand{\corr}{{\rm corr}}

\newcommand{\bfTh}{\mbox{\boldmath$\Theta$}}

\newcommand{\rhs}{right-hand side}
\newcommand{\df}{distribution function}

\newcommand{\beao}{\begin{eqnarray*}}
\newcommand{\eeao}{\end{eqnarray*}\noindent}

\newcommand{\beam}{\begin{eqnarray}}
\newcommand{\eeam}{\end{eqnarray}\noindent}

\newcommand{\beqq}{\begin{equation}}
\newcommand{\eeqq}{\end{equation}\noindent}

\newcommand{\bce}{\begin{center}}
\newcommand{\ece}{\end{center}}

\newcommand{\barr}{\begin{array}}
\newcommand{\earr}{\end{array}}

\newcommand{\stv}{\stackrel{v}{\rightarrow}}
\newcommand{\stw}{\stackrel{w}{\rightarrow}}

\newcommand{\vague}{\stackrel{\lower0.2ex\hbox{$\scriptscriptstyle
                    \it{v} $}}{\rightarrow}}
\newcommand{\weak}{\stackrel{\lower0.2ex\hbox{$\scriptscriptstyle
                    \it{w} $}}{\rightarrow}}
\newcommand{\what}{\stackrel{\lower0.2ex\hbox{$\scriptscriptstyle
                    \it{\hat{w}} $}}{\rightarrow}}

\newcommand{\bdis}{\begin{displaymath}}
\newcommand{\edis}{\end{displaymath}\noindent}

\newcommand{\nto}{n\to\infty}

\newcommand{\xto}{x\to\infty}

\newcommand{\wh}{\widehat}

\newcommand{\la}{\lambda}

\newcommand{\regvary}{regularly varying}
\newcommand{\slvary}{slowly varying}
\newcommand{\regvar}{regular variation}

\newcommand{\bbr}{{\mathbb R}}

\newcommand{\bbz}{{\mathbb Z}}

\newcommand{\bbs}{{\mathbb S}}

\newcommand{\con}{convergence}

\newcommand{\st}{such that}
\newcommand{\fif}{if and only if}
\newcommand{\wrt}{with respect to}

\newcommand{\fct}{function}

\newcommand{\ds}{distribution}

\newcommand{\rep}{representation}

\newcommand{\seq}{sequence}

\newcommand{\ms}{measure}

\newcommand{\bfx}{{\bf x}}
\newcommand{\bfX}{{\bf X}}
\newcommand{\bfB}{{\bf B}}
\newcommand{\bfY}{{\bf Y}}

\newcommand{\bfA}{{\bf A}}

\newcommand{\bfZ}{{\bf Z}}
\newcommand{\bfW}{{\bf W}}

\newcommand{\bfa}{{\bf a}}

\newcommand{\bfs}{{\bf s}}

\newcommand{\bali}{\begin{align}}
\newcommand{\eali}{\end{align}}

\begin{document}

\bibliographystyle{alpha}
\title[The extremogram and the cross-extremogram of multivariate high
frequency data]
      {The extremogram and the cross-extremogram for a bivariate
        \garch\ process}
\today
\author[M. Matsui]{Muneya Matsui}
\author[T. Mikosch]{Thomas Mikosch}
\address{Department of Business Administration, Nanzan University,
18 Yamazato-cho Showa-ku Nagoya, 466-8673, Japan}
\email{mmuneya@nanzan-u.ac.jp}
\address{T. Mikosch, University of Copenhagen, Department of Mathematics,
Universitetsparken 5,
DK-2100 Copenhagen\\ Denmark} \email{mikosch@math.ku.dk}

\begin{abstract} 
In this paper, we derive some \asy\ theory for the extremogram and
cross-extremogram of a bivariate \garch\ process. We show that the
tails of the components of a bivariate \garch\ process may exhibit 
power law behavior but, depending on the choice of the parameters, 
the tail indices of the components may differ. We apply the theory to
5-minute return data of stock prices and foreign exchange rates. 
We judge the fit of a bivariate \garch\ model by considering the
sample extremogram and cross-extremogram of the 
residuals. The results are in agreement with the iid hypothesis of the
two-dimensional innovations \seq . The cross-extremograms at lag zero
have a value significantly distinct from zero. This fact points at 
some strong extremal dependence of the components of the innovations.
\end{abstract}

\thanks{
Muneya Matsui's research is partly supported by the JSPS Grant-in-Aid
for Young Scientists B (25870879).
}
\maketitle 
\section{The extremogram and the cross-extremogram}
In this paper we conduct an empirical study of extremal serial dependence
in a bivariate return series. Our main tools for describing extremal dependence
will be the {\em extremogram} and the {\em cross-extremogram}. For the
sake of argument and for
simplicity, we restrict ourselves to bivariate series $\bfX_t=
(X_{1,t},X_{2,t})'$, $t\in \bbz$. We assume that $(\bfX_t)$ has the
following structure:
\beam\label{eq:13}
\bfX_t= \Sigma_t\, \bfZ_t\,,\qquad t\in\bbz\,, 
\eeam
where $(\bfZ_t)$ constitutes an iid bivariate noise \seq\ and
\beao
\Sigma_t={\rm diag}(\sigma_{1,t},\sigma_{2,t})\,,\quad t\in\bbz\,, 
\eeao
where $\sigma_{i,t}$ is the (non-negative) volatility of $X_{i,t}$.
We will assume that $(\bfX_t)$, $(\Sigma_t)$ constitute strictly
stationary \seq s and that $\Sigma_t$ is predictable \wrt\ the
filtration generated by $(\bfZ_s)_{s\le t}$. We also assume that
$\bfZ_t=(Z_{1,t},Z_{2,t})'$ has mean zero and 
covariance matrix  (standardized to correlations) 
\beam\label{eq:14}
P =\left(\barr{ll}
1&\rho \\
\rho& 1
\earr\right)\,,
\eeam
where $\rho=\corr(Z_{1,t},Z_{2,t})$. Later, we will choose parametric
models for $(\bfX_t)$ such as univariate \garch\ models both for $X_{i,t}$,
$i=1,2,$ or a vector GARCH(1,1) model; see Section~\ref{subsec:garch} 
for model descriptions. In the context of these parametric models, the
choice of the covariance matrix $P$ as a correlation matrix is a
matter of identifiability of the model since one can always swap a
positive constant multiplier between $\Sigma_t$ and~$\bfZ_t$.
\par
The extremogram and cross-extremogram for a stationary \seq\
$(\bfX_t)$ were introduced in Davis and Mikosch
\cite{davis:mikosch:2009} and further developed in 
Davis et
al. \cite{davis:mikosch:cribben:2012,davis:mikosch:zhao:2013}.
In standardized forms, the extremogram and cross-extremogram 
of a bivariate \seq\ $(\bfX_t)$ are given by 
\beao
\left(\barr{ll}\rho_{11}(h)& \rho_{12}(h)\\
\rho_{21}(h)& \rho_{22}(h)\earr\right)\,,\qquad h=0,1,2,\ldots,
\eeao
where 
\beam\label{eq:2}
\left(\barr{l}
\rho_{11}(h)\\
\rho_{22}(h)\\
\rho_{12}(h)\\
\rho_{21}(h)
\earr\right)=
 \lim_{\xto}\left(\barr{l} \P( X_{1,h}\in x\,A\mid X_{1,0}\in x \,
A)\\
\P( X_{2,h}\in x\,B\mid X_{2,0}\in x \,
B)\\
\P( X_{2,h}\in x\,B\mid X_{1,0}\in x \,
A)\\
\P( X_{1,h}\in x\,A\mid X_{2,0}\in x \,
B)\earr\right). 
\eeam
Here $A,B$ are sets bounded away from zero. Typically,  we choose intervals
$(1,\infty)$, $(-\infty,-1)$ for $A,B$ and we also suppress the
dependence on $A,B$ in the $\rho_{ij}$-notation. Notice that $x
A=\{xy: y\in A\}$ has interpretation as an extreme event if $x$ is
sufficiently large. Thus, the extremogram $\rho_{ii}(h)$ describes the
likelihood of an extreme event at lag $h$ given there is an extreme
event in the $i$th component at time zero. Correspondingly, the 
cross-extremogram $\rho_{ij}(h)$ for $i\ne j$ describes the likelihood
of an extreme event at  time lag $h$ in the $j$th component given there
is an extreme event at time zero in the $i$th component. In general,
$\rho_{12}(h)\ne \rho_{21}(h)$.
The limits
$\rho_{ij}$ can be understood as generalizations of the
(upper) 
tail dependence
coefficient 
$\rho=\lim_{q\uparrow 1}\P(Y>F_Y^{\leftarrow}(q)\mid
X>F_X^{\leftarrow}(q))$ for a bivariate vector $(X,Y)$ to the \ts\ context.
Here $F_X^{\leftarrow}(q),F_Y^{\leftarrow}(q)$  are the $q$-quantiles
of the \ds s of $X,Y$, respectively. The tail dependence coefficients
have been propagated for measuring extremal dependence in a bivariate
vector in the context of quantitative risk management; see e.g.
McNeil et al. \cite{mcneil:frey:embrechts:2005}. 
\par
Moreover, each of the quantities $\rho_{ij}(h)$ has interpretation as
a limiting covariance or cross-covariance \fct . For example.
\beao
\rho_{11}(h)&=& \lim_{\xto} \dfrac{\cov(\1(X_{1,0}\in x\,A), \1(X_{1,h}\in x\,A))+
[\P(X_{1,0}\in x\,A)]^2}{\P(X_{1,0}\in  x\,A)}\\
&=&
\lim_{\xto} \dfrac{\cov(\1(X_{1,0}\in x\,A), \1(X_{1,h}\in x\,A))}{\P(X_{1,0}\in  x\,A)}
\\
&=&
 \lim_{\xto} \P(X_{1,h}\in x\,A\mid X_{1,0}\in x\,A) \,.
\eeao
The limits $\rho_{ij}(h)$ in \eqref{eq:2} do not automatically exist. A
convenient theoretical assumption for their existence is the condition
of {\em \regvar\ of the \ts\ $(\bfX_t)$}. This notion 
is explained in Section~\ref{subsec:regvar}. Its close relationship with
GARCH models is investigated in Section~\ref{subsec:garch}.
 Return series are often heavy-tailed and
therefore it is attractive to model them by a \regvary\ model.
For example, under mild conditions the GARCH model automatically 
ensures that sufficiently high moments
of the series are infinite. In  particular, univariate and
multivariate GARCH models exhibit power laws. This will be
explained in Section\ref{subsec:garch}.
 
\section{Some preliminaries}
\subsection{Regularly varying \ts }\label{subsec:regvar}
We say that an $\bbr^d$-valued strictly 
stationary \ts\ $(\bfX_t)$ is {\em \regvary\ 
with index $\alpha>0$} if its \fidi s are \regvary\ in the following
sense:
for every $h\ge 0$, the following limits in \ds\ exist
\beam\label{eq:jan6a}
\P(x^{-1} (\bfX_0,\ldots,\bfX_h)\in \cdot\mid |\bfX_0|>x )\stw
\P( (\bfY_0,\ldots,\bfY_h)\in \cdot)\,,\qquad \xto\,,
\eeam 
where the limit vector $(\bfY_0,\ldots,\bfY_h)$ has the same \ds\ as 
$|\bfY_0|(\bfTh_0,\ldots,\bfTh_h)$, the \ds\ of $|\bfY_0|$ is given by
$\P(|\bfY_0|>y)=y^{-\alpha}$, $y>1$, and $|\bfY_0|$ and 
$(\bfTh_0,\ldots,\bfTh_h)$ are independent. Of course,  the \ds\ of 
$\bfTh_0$, the {\em spectral \ms }, is concentrated on the unit sphere
$\bbs^{d-1}=\{\bfx \in\bbr^d: |\bfx|=1\}$. The spectral \ms\
describes the likelihood of the directions of extreme values of the
lagged vector $\bfX_0$. Here $\stw$ denotes weak convergence and $|\cdot|$ denotes any norm in $\bbr^d$;
from now on we choose the Euclidean one. The aforementioned
definition of a \regvary\ \ts\ is based on work by Basrak and Segers
\cite{basrak:segers:2009} which yields a convenient description of
the topic.  
Davis and Hsing \cite{davis:hsing:1995} introduced the 
notion of a \regvary\ \ts\ which is attractive for describing
serial extremal dependence in the presence of heavy tails. They 
used an alternative definition of multivariate \regvar\ which is 
equivalent to the definition above.
\par
A direct con\seq\  of the \regvar\ of a \ts\ is that 
\beam\label{eq:3}
\mbox{$\P(|\bfX_0|>x)=x^{-\alpha}L(x)$, \qquad $x>0$, \qquad for a \slvary\ \fct\
  $L$,}
\eeam
i.e., $L$ is a positive \fct\ on $(0,\infty)$ \st\ 
$L(cx)/L(x)\to 1$  as $\xto$ for any $c>0$. 
Then we also have
\beam\label{eq:4}
\P(\bfX_0/|\bfX_0|\in \cdot\mid |\bfX_0|>x)\stw
\P(\bfTh_0\in\cdot)\,,\qquad \xto\,.
\eeam
Regular variation of the marginal \ds\ of the \ts\ is equivalent
to the set of relations \eqref{eq:3} and \eqref{eq:4}.
A further con\seq\ is
that $\P(\bfs '\bfX_0>x)/\P(|\bfX_0|>x)\to e_\alpha(\bfs)$ as $\xto$ for any
choice of $\bfs\in \bbs^{d-1}$ and some \fct\ 
$e_\alpha$ \st\ $e_\alpha(\bfs_0)\ne 0$ for some $\bfs_0\in\bbs^{d-1}$ and $e_\alpha(t\bfs)=t^{-\alpha}e_\alpha(\bfs)$, $t>0$. For proofs 
of the aforementioned properties and further reading on \regvar , we refer to 
Resnick \cite{resnick:1987,resnick:2007}.
\par
A particular con\seq\ of the property of \regvar\ of a \ts\ $(\bfX_t)$
is the fact that the limits in  \eqref{eq:2}, leading to
the extremogram and cross-extremogram, are well defined.
For this reason, we will assume that $(\bfX_t)$ is \regvary\ or we will
assume that a deterministic monotone 
increasing transformation of the components $X_{i,t}$, $i=1,2$, of
$\bfX_t$ results in a \regvary\ \ts . Such transformations can be
necessary, for example, if both components are not \regvary\ or if both
components have rather different tail behavior. These cases are 
relevant for real-life \ts . For the sake of argument,
assume that $(\bfX_t)$ is a bivariate strictly stationary Gaussian \ts .
This is not a \regvary\ time series. However, the extremogram and
cross-extremogram of this \seq\ exist for various sets $A,B$, for
example, if $A=B=(1,\infty)$ (a corresponding remark applies if 
$A$ or $B$ is the set $(-\infty,-1)$). If $G$ denotes the \df\ of  
a $t$-\ds\ with $\alpha$ degrees of freedom then 
calculation yields that 
\beam\label{eq:5}
\big(
G(F_{X_{1,0}}(X_{1,t})),
G(F_{X_{2,0}}(X_{2,t}))\big)\,,\qquad t\in\bbz\,,
\eeam
has $G$-distributed marginals and  one can indeed show that the
transformed \ts\ is  \regvary\ with index $\alpha$. The same
transformation arguments apply to a non-Gaussian  \ts .
In contrast to a Gaussian \ts , in general one cannot ensure that the
resulting \ts\ is \regvary\ in the sense defined above. 
Given that a
transformation of the type  \eqref{eq:5} yields a \regvary\ \ts , one
can modify the cross-extremogram e.g. for the sets $A=B=(1,\infty)$ in
the following way:
\beao
\rho_{12}(h)&=& \lim_{q\uparrow 1}\P( X_{2,h} > F_{X_{2,0}}^\leftarrow
(q)\mid  X_{1,0} > F_{X_{1,0}}^\leftarrow
(q))\\
&=& \lim_{x\to\infty }\P( G(F_{X_{2,0}}(X_{2,h})) > x\mid  
G(F_{X_{1,0}}(X_{1,0})) > x)\,.
\eeao
For practical purposes, we will replace the high quantiles 
$F_{X_{i,0}}^\leftarrow (q)$, $i=1,2$, by their empirical versions,
such as the $97\%$, $98\%$,... componentwise empirical quantiles, 
depending on the
sample size available.
\par
Regular
variation of a \ts\ is a convenient theoretical property but it 
cannot be tested on data. Therefore we will assume  a GARCH model for
$(\bfX_t)$. This model assumption ensures \regvar .

\subsection{Univariate \garch\ models}\label{subsec:garch}
From Bollerslev \cite{bollerslev:1986} 
recall the notion of a univariate \garch\ model 
\beam\label{eq:9}
X_t= \sigma_t\,Z_t\,, \qquad t\in\bbr\,,
\eeam
where $(Z_t)$ is an iid unit variance mean zero \seq\ and $(\sigma_t)$
is a positive volatility \seq\ whose dynamics is given by the causal
non-zero solution to the \sre
\beam\label{eq:6}
\sigma_t^2 =\alpha_0+
\alpha_1\,X_{t-1}^2+\beta_1\sigma_{t-1}^2= 
\alpha_0+(\alpha_1Z_{t-1}^2+\beta_1)\,\sigma_{t-1}^2\,,\qquad t\in\bbz\,.
\eeam
Here $\alpha_0>0$ and $\alpha_1>0$, $\beta_1\ge 0$ are constants. 
The probabilistic structure of $(\sigma_t^2)$ can be investigated in
the context of  solutions to the general \sre\ 
\beam\label{eq:8}
Y_t= A_t\,Y_{t-1}+B_t\,,\qquad t\in\bbr\,,
\eeam
where $(A_t,B_t)$, $t\in\bbz$, constitutes an $\bbr_+^2$-valued iid \seq
.
Indeed, $(\sigma_t^2)$ satisfies this equation with $B_t=\alpha_0$
and $A_t=\alpha_1Z_{t-1}^2+\beta_1$. Based on the theory for these equations
(see 
Bougerol and 
Picard \cite{bougerol:picard:1992}), we conclude that a strictly stationary 
positive solution $(\sigma_t^2)$ to \eqref{eq:6} exists \fif\ 
\beam\label{eq:7}
\E \log (\alpha_1Z_0^2+\beta_1)<0\qquad\mbox{and}\qquad \alpha_0>0\,.
\eeam In view of
Jensen's inequality and since $\E Z_0^2=1$,  
$\E \log (\alpha_1Z_0^2+\beta_1)\le\log  \E
(\alpha_1 Z_0^2+\beta_1)= \log(\alpha_1+\beta_1)$. Therefore the condition
$\alpha_1+\beta_1<1$ ensures strict stationarity as well as second
order stationarity of $(\sigma_t)$ and $(X_t)$, 
but the condition  \eqref{eq:7} is much more general and also allows
for certain choices of $\alpha_1,\beta_1$ \st\ $\alpha_1+\beta_1\ge
1$; see Nelson \cite{nelson:1990}, Bougerol and 
Picard \cite{bougerol:picard:1992}. In the latter cases, $ \E[X_0^2]=\infty$. 
\par
The solution to \eqref{eq:8} has a rather surprising property which
was discovered by Kesten \cite{kesten:1973}; see also Goldie
\cite{goldie:1991}. Under mild conditions, there exists a positive
constant $c_0$ \st\ $\P(Y_0>x)\sim c_0\,x^{-\alpha}$ for some 
$\alpha>0$. We apply the aforementioned theory to \eqref{eq:6} and 
get the following result which can be found in Goldie's
\cite{goldie:1991} paper as regards the marginal \ds s.
Mikosch and \sta\ \cite{mikosch:starica:2000} proved that
$(\sigma_t)$ and $(X_t)$ are \regvary\ \ts .
\bpr\label{prop:1}
Assume that $\alpha_0>0$, $Z_0$ has Lebesgue density and there exists $\alpha>0$
\st\ 
\beam\label{eq:12}
\E (\alpha_1Z_0^2+\beta_1)^{\alpha/2}=1\,,
\eeam 
and 
$\E \big[(\alpha_1Z_0^2+\beta_1)^{\alpha/2}
\log^+(\alpha_1Z_0^2+\beta_1)\big]<\infty$. Then there exist a unique
strictly stationary causal non-zero solutions to \eqref{eq:6} and \eqref{eq:9},
and there exists a
constant
$c_0>0$ \st\ 
\beam\label{eq:11}
\P(\sigma_0>x)\sim c_0\,x^{-\alpha}\,,\qquad\xto\,.
\eeam
Moreover, as $\xto$,
\beam\label{eq:10}
\P(X_0>x)\sim \E [Z_0^+]\,\P(\sigma_0>x)\quad\mbox{and}\quad 
\P(X_0\le -x)\sim \E [Z_0^-]\,\P(\sigma_0>x)\,,
\eeam
where $x^\pm=\max(0,\pm x)$. In addition, the \seq s $(\sigma_t)$ and
$(X_t)$ are \regvary\ with index $\alpha$.
\epr
Relation \eqref{eq:10} is an immediate con\seq\ of \eqref{eq:11} and 
a result of Breiman \cite{breiman:1965} about the tails of products of
independent \rv s; cf. Jessen and Mikosch \cite{jessen:mikosch:2006}.
\section{Bivariate \garch\ processes and their properties}
Our next goal is to consider multivariate extensions of the \garch\
model of the type described in \eqref{eq:13}. A simple way of doing
this is by assuming that both component \seq s  $(X_{i,t})$, $i=1,2$, 
constitute univariate \garch\ processes, i.e., $(\bfX_t)$ in
\eqref{eq:13} is specified via the vector recursion
\beao
\left(\barr{l}\sigma_{1,t}^2\\
\sigma_{2,t}^2\earr\right)&=&
\left(\barr{l}\alpha_{01}\\\alpha_{02}\earr\right)
+\left(\barr{cc}
\alpha_{11}& 0\\
0&\alpha_{22}\earr\right)\,\left(\barr{l}X_{1,t-1}^2\\
X_{2,t-1}^2\earr\right)+\left(\barr{cc}
\beta_{11}& 0\\
0&\beta_{22}\earr\right)\,\left(\barr{l}\sigma_{1,t-1}^2\\
\sigma_{2,t-1}^2\earr\right)\\
&=&
\left(\barr{l}\alpha_{01}\\\alpha_{02}\earr\right)+
\left(\barr{cc}\alpha_{11}Z_{1,t-1}^2+\beta_{11}&0\\
0&\alpha_{22}Z_{2,t-1}^2+\beta_{22}
\earr\right)\,\left(\barr{l}\sigma_{1,t-1}^2\\
\sigma_{2,t-1}^2\earr\right)\,,
\eeao 
and $(\bfZ_t)$ is an iid \seq\ with covariance matrix 
$P$ given in \eqref{eq:14}.
We can apply the univariate theory to the components  $(\sigma_{i,t}^2)$, $i=1,2$.
There exist unique strictly stationary solutions $(\sigma_{i,t}^2)$, $i=1,2$,
\fif\ $\alpha_{0i}\ne 0$ and 
$\E \log^+(\alpha_{ii}Z_{i,0}^2+\beta_{ii})<0$ for $i=1,2$.
Then we may also conclude that
the resulting unique bivariate processes $(\Sigma_t)$ and $(\bfX_t)$ are 
strictly stationary. Notice that the dependence structure 
between the univariate component
processes is then completely determined by the dependence 
structure of the components of the noise $(\bfZ_t)$. 
We can also apply
Proposition~\ref{prop:1} to get conditions for power law tails and
\regvar\ of the component processes of $(\bfX_t)$. 
\bre\label{rem:1}
The crucial condition for the componentwise
tail behavior is \eqref{eq:12}. Since the \ds s of $Z_{i,0}$, $i=1,2$,
and the parameter sets $(\alpha_{ii},\beta_{ii})$, $i=1,2$, may be distinct,
$X_{1,t}$ and $X_{2,t}$ will in general have different tail indices
$\alpha_1$ and $\alpha_2$, respectively. This fact can be considered
an advantage when studying multivariate return series. Indeed,
it is not likely that the tail indices of the univariate components of real-life \ts\ coincide.

\ere
\par
There exist various extensions of a univariate GARCH
model to the multivariate case. We stick here to the 
{\em constant conditional correlation model} of Bollerslev
\cite{bollerslev:1990}.
It is the model \eqref{eq:13} with specification
\beam
\begin{split}\label{eq:15}
\left(\barr{l}\sigma^2_{1,t}  \\  
\sigma^2_{2,t}\earr
\right)
&=& \left(
\barr{l}\a_{01}  \\\a_{02}   \earr\right)
+\left(\barr{cc}\a_{11} & \a_{12}  \\
      \a_{21} & \a_{22}\earr \right)\, 
\left(\barr{l}X_{1,t-1}^2  \\X_{2,t-1}^2   \earr\right)
 + \left(\barr{cc}\b_{11} & \b_{12}  \\\b_{21} & \b_{22} \earr
 \right)\,\left(\barr{c}\sigma^2_{1,t-1}  \\\sigma^2_{2,t-1}\earr
  \right)\\
&=& \left(
\barr{l}\a_{01}  \\\a_{02}   \earr\right)+\left(\barr{cc}\alpha_{11}Z_{1,t-1}^2+\beta_{11}&\alpha_{12}Z_{2,t-1}^2+
\beta_{12}\\
\alpha_{21}Z_{1,t-1}^2+\beta_{21}& \alpha_{22}Z_{2,t-1}^2+\beta_{22}
\earr\right)\,\left(\barr{l}\sigma_{1,t-1}^2\\\sigma_{2,t-1}^2\earr
\right)\,.
\end{split}
\eeam
Writing $\bfW_t=(\sigma^2_{1,t} \,,
\sigma^2_{2,t})'$, 
\beam\label{eq:20a}
\bfB_t= \left(
\barr{l}\a_{01}  \\\a_{02}   \earr\right)\quad\mbox{and}\quad
\bfA_t=\left(\barr{cc}\alpha_{11}Z_{1,t-1}^2+\beta_{11}&\alpha_{12}Z_{2,t-1}^2+
\beta_{12}\\
\alpha_{21}Z_{1,t-1}^2+\beta_{21}& \alpha_{22}Z_{2,t-1}^2+\beta_{22}
\earr\right)\,,
\eeam
we see that we are again in the framework of the \sre\ \eqref{eq:8},
but this time for vector-valued $\bfB_t$ and matrix-valued $\bfA_t$:
\beam\label{eq:jan6b}
\bfW_t=\bfA_t\,\bfW_{t-1}+\bfB_t\,,\qquad t\in\bbz\,.
\eeam
Kesten \cite{kesten:1973} also provided the corresponding theory  
for stationarity and tails in this case. \sta\ \cite{starica:1999}
dealt with the corresponding problems for vector \garch\ processes,
making use of the theory in Kesten \cite{kesten:1973},
Bougerol and Picard \cite{bougerol:picard:1992}
and its
specification to the tails of GARCH models 
in Basrak et al.~\cite{basrak:davis:mikosch:2002}. In the bivariate 
\garch\ case the theory in \sta\ \cite{starica:1999} can be written in a
more compact form due to the \rep\ \eqref{eq:15}; in the case of
higher order GARCH models    \eqref{eq:15} has to be written as an
equation for vectors involving both $\sigma^2$- and $X^2$-terms at
more than 1 lag.
\par
According to Bougerol and Picard \cite{bougerol:picard:1992}, 
 \eqref{eq:15} has a unique strictly stationary solution if
the {\em top Lyapunov exponent} $\gamma$ 
associated with the \seq\ $(\bfA_t)$ is negative, i.e.,
\beam\label{eq:16}
\gamma=\lim_{\nto} n^{-1} \log \|\bfA_1\cdots \bfA_n\|<0\,,
\eeam
where $\|\cdot\|$ denotes the matrix norm and
the limit on the \rhs\ exists a.s. In view of Remark on p.~122
in \cite{bougerol:picard:1992},
a sufficient
condition for $\gamma<0$ is that the matrix 
\beam\label{eq:spectral}
\E \bfA_1= \left(\barr{cc}\a_{11}+\b_{11}& \a_{12}+\b_{12}\\
\a_{21}+\b_{21}& \a_{22}+\b_{22}\earr\right) =:
\left(\barr{cc}a_{11} & a_{12}  \\
      a_{21} & a_{22}\earr \right) 
\eeam
has spectral radius smaller than 1. We assume that all entries 
of $\E \bfA_1$ are positive.  Then,
by the Perron--Frobenius theorem (see Lancaster \cite{lancaster:1969}, Section 9.2),  $\E \bfA_1$ has a dominant single positive eigenvalue. Keeping this fact 
in mind, the largest positive solution to the characteristic equation 
${\rm det}\,(\la I-\E \bfA_1)=0$ yields 
the sufficient condition 
\beam\label{eq:stationarty}
&&\dfrac{a_{11}+a_{22}}{2}+ 
\sqrt{\Big(\dfrac{a_{11}-a_{22}}{2}\Big)^2+ a_{12}a_{21}}\quad <\,1\,.
\eeam
\par
Next we give sufficient conditions for the \regvar\ of a bivariate \garch\
process $(\bfX_t)$. The proof is based on Kesten's fundamental results
\cite{kesten:1973}, in particular Theorem~4. \sta\ \cite{starica:1999} gave a similar result,
referring to Basrak et al.~\cite{basrak:davis:mikosch:2002} for a
related proof in the situation of a univariate \garchpq\ process. 
We give a proof by verifying Kesten's conditions.
\bpr\label{prop:2} Consider the bivariate \garch\ model and assume the following conditions:
\begin{enumerate}
\item
Condition \eqref{eq:16}. 
\item
$\bfZ_0$ has Lebesgue density in $\bbr^2$.
\item
There exists $p>0$ \st\
\beam\label{eq:17}
\E [|\bfZ_0|^{2p}\log ^+ |\bfZ_0|]<\infty\quad\mbox{and}\quad
\E \big[\min_{i=1,2}\big(\sum_{j=1}^2(\alpha_{ij}Z_{j,0}^2+\b_{ij})\big)^p\big]
\ge 2^{p/2}\,.\eeam
\item
All entries of $\bfA_0$ are positive a.s., 
$\alpha_{0i}>0$, $i=1,2$, and not all values $\alpha_{ij}$, $1\le i,j\le 2$, vanish. 
\end{enumerate}
Then there exists a unique $\alpha\in (0,2\,p]$ \st\ 
\beam
0= \lim_{\nto} n^{-1}\,\log \E \big[\|\bfA_1\cdots \bfA_n\|^{\alpha/2}\big]\,,
\label{eq:jan6c}\eeam
there exists a strictly stationary causal non-zero solution $(\bfX_t)$ to
\eqref{eq:13} with specification \eqref{eq:15} and $(\bfX_t)$ is
\regvary\ with index $\alpha$. In particular, for every $n\ge 1$, there exists a non-null Radon \ms\ 
$\mu$ on $\bbr^{2n}\backslash\{{\bf0}\}$ \st
\beao
x^\alpha \,\P\big(x^{-1}(\bfX_1,\ldots,\bfX_n)\in\cdot)\stv  \mu_n(\cdot)\,,\qquad \xto\,.
\eeao
Here $\stv$ denotes vague \con .
\epr
\begin{proof}
According to Kesten \cite{kesten:1973}, Theorem 4, there exist
\begin{itemize}
\item
a
unique strictly stationary solution $(\bfW_t)$ to the equation \eqref{eq:jan6b},
\item
a positive value $\alpha$ and a non-negative \fct\ $e_\alpha$ on $\bbs^1$ \st
\beam\label{eq:18}
\lim_{\xto}x^{\alpha/2}\P(\bfu'\bfW_0>x)= e_\alpha(\bfu)\,,\qquad \bfu\in\bbs^1\,,
\eeam
and the \fct\ $e_\alpha$ is positive for $\bfu\in \bbs^1$ \st\ 
$\bfu\ge {\bf0}$,\end{itemize}
if the following conditions hold:
\begin{enumerate}
\item[\rm 1.]
$\bfA_0\ge {\bf0}$ and $\bfB_0\ge {\bf0}$ and $\bfB_0\ne {\bf0}$, where
	     $\bfC\ge {\bf0}$ (resp. $>{\bf0}$) implies all entries in
	     $\bfC$ are non-negative (resp. positive).
\item[\rm 2.]
The additive group generated by the numbers $\log \rho(\bfa_1\cdots \bfa_n)$
is dense in $\bbr$, where $\bfa_i$ are elements in the support of the
\ds\ of $\bfA_0$ \st\ $\bfa_1\cdots \bfa_n$ has positive entries 
 and $\rho$ is the spectral radius.
\item[\rm 3.] Condition  \eqref{eq:16} holds.
\item[\rm 4.] There is $\alpha>0$ \st\ \eqref{eq:jan6c} holds.
\item[\rm 5.] $\E[\|\bfA_0\|^{\alpha/2}\log^+\|\bfA_0\|]<\infty$ and $\E[|\bfB_0|^{\alpha/2}]<\infty$.

\end{enumerate}
Condition 1 holds in view of the assumptions $\bfA_0>{\bf0}$ a.s. and $\bfB_0>{\bf0}$.\\
Condition 2. We assume that $\bfA_0>{\bf0}$ a.s. 
Therefore $\bfa_1\cdots\bfa_n>{\bf0}$ for any $n\ge 1$ and any $\bfa_i$ in the support of $\bfA_0$. Since we assume a Lebesgue
density for $\bfZ_0$ there exists an open set in $\bbr^2$ where this density is 
positive. Therefore 
and since not all values $\alpha_{ij}$ vanish,
there exists a continuum of values $\rho(\bfa_1)$ for $\bfa_1$ in
the support of $\bfA_0$.\\
Conditions 3 and 5 follow from the assumptions.\\
Condition 4. The existence of such an $\alpha$ 
 follows from the existence of $p>0$ \st\ \eqref{eq:17} holds. Then $\alpha\le p$.
\par
Thus Kesten's Theorem 4 can be applied. In particular, \eqref{eq:18} holds.
Due to results in
Boman and Lindskog \cite{boman:lindskog:2009} and since $\bfW_0$ is positive, 
\eqref{eq:18} implies that $\bfW_0$ is \regvary\ with index $\alpha$ in
the sense of  Section~\ref{subsec:regvar}. 
\par
Next we show that the
\fidi s of $(\bfW_t)$ 
are \regvary . This follows from the \rep\
\beam\label{eq:20}
(\bfW_1,\ldots,\bfW_t)=({\mathbf \Pi}_1,\ldots,{\mathbf \Pi}_t)\,\bfW_0+ (\bfR_1,\ldots,\bfR_t)\,,\qquad t\ge 1\,,
\eeam 
where ${\mathbf \Pi}_t=\bfA_t\cdots \bfA_1$ for $t\ge 1$,
$({\mathbf \Pi}_1,\ldots,{\mathbf \Pi}_t)$, $(\bfR_1,\ldots,\bfR_t)$ have moment of order $\alpha/2$ and 
are independent of $\bfW_0$. Now
an application of the multivariate Breiman theorem in Basrak et al.
\cite{basrak:davis:mikosch:2002a} yields the \regvar\ of the \fidi s of
$(\bfW_t)$ with index $\alpha/2$, due to the \regvar\ of
$\bfW_0$ with the same index. Hence $(\Sigma_t)=(({\rm diag}(\bfW_t))^{1/2})$ inherits \regvar\ with index $\alpha$. Here $\bfx^{1/2}$ for any matrix 
or vector $\bfx$ refers to taking square roots for all entries. 
\par
It remains to show that $(\bfX_t)$ is \regvary\ with index $\alpha$.
We write
$\wt \Sigma_t=({\rm diag} (\bfW_t-\bfR_t))^{1/2}$. It is not difficult to see that 
$|(\Sigma_t-\wt\Sigma_t)\bfZ_t|$ is dominated by $c|\bfR_t|^{1/2} \,|\bfZ_t|$ 
for some constant $c$ and this
bound has finite $\alpha$th moment. Therefore
\beao
\lim_{\xto}x^{\alpha}\,\P \big( \big|(\Sigma_1 \bfZ_1,\ldots,\Sigma_t\bfZ_t)-
(\wt \Sigma_1 \bfZ_1,\ldots,\wt \Sigma_t\bfZ_t)\big|>x\big)=0\,,\qquad \xto\,.
\eeao
Since $\bfW_0$ is \regvary\ with index $\alpha/2$ an application of the 
multivariate Breiman result (see Basrak et al. 
\cite{basrak:davis:mikosch:2002a}) shows that $ ({\mathbf \Pi}_1,\ldots,{\mathbf \Pi}_t)\bfW_0$ is \regvary\
with index $\alpha/2$ as well.
Combining these facts, we conclude that
\beao
(\wt \Sigma_1 \bfZ_1,\ldots,\wt \Sigma_t\bfZ_t)\quad\mbox{and}\quad 
(\Sigma_1 \bfZ_1,\ldots, \Sigma_t\bfZ_t)
\eeao
have the same tail behavior and are \regvary\ with index $\alpha$; cf. 
Jessen and Mikosch \cite{jessen:mikosch:2006}.  In particular, we have
\beao
\lefteqn{\P\big(x^{-1}\big(({\rm diag} ({\mathbf \Pi}_1\bfW_0))^{1/2} \bfZ_1,\ldots,
({\rm diag} ({\mathbf \Pi}_t\bfW_0))^{1/2}\bfZ_t\big)\in\cdot\mid |\bfW_0|>x\big)}\\
&\stw&
 \P(Y_0\,\big(({\rm diag} ({\mathbf \Pi}_1\bfTh_0))^{1/2} \bfZ_1,\ldots,
({\rm diag} ({\mathbf \Pi}_t\bfTh_0))^{1/2}\bfZ_t\big)\in\cdot\big)\,,
\eeao
where $\P(Y_0>x)=x^{-\alpha}$, $x>1$, $Y_0$ is independent of
 $\bfTh_0,\bfZ_1,\ldots,\bfZ_t$ and $\bfTh_0$ has the spectral \ds\ of
 $\bfW_0$. 
\end{proof}
\bre
In view of Kesten's result, relation \eqref{eq:18} holds for any 
$\bfu\in\bbs^1$ and $e_\alpha(\bfu)\ne 0$ for $\bfu\ge {\bf0}$. 
In particular, for $\bfu_1=(0,1)$ and $\bfu_2=(1,0)$ we
conclude that $\P(\sigma_{i,0}>x)\sim c_i\,x^{-\alpha}$ as $\xto$, 
where both constants $c_i$ are positive. In turn, 
Breiman's result \cite{breiman:1965} ensures that 
\beao
\P(X_{i,0}^\pm >x)\sim \E
[(Z_{i,0}^\pm)^\alpha]\,\P(\sigma_{i,0}>x)\,,\qquad \xto\,,\qquad i=1,2.
\eeao
This means that the right and left tails of the \ds\ of 
$\bfX_0$ are equivalent and they have the same tail index
$\alpha$. This is in contrast to the case when
$\alpha_{ij}=\beta_{ij}=0$ for $i\ne j$ (see Remark~\ref{rem:1}),
where the components of $\bfX_t$ may have different tail behavior.
From a modeling point of view, this property does not allow for much
flexibility as regards the componentwise extremes in a multivariate
return series. This fact can be understood as a disadvantage
of the constant conditional correlation model \wrt\ more realistic
modeling of the extremes of multivariate return models.
\par
The crucial condition in Proposition~\ref{prop:2} which makes the
difference to Proposition~\ref{prop:1} is the assumption that all
entries of $\bfA_0$ must be positive and random. 
This condition is also satisfied if 
$\a_{ii}>0$, $i=1,2,$ $\alpha_{ij}=0$ and $\beta_{ij}>0$, $i\ne j$,
i.e., the  off-diagonal elements in the matrix $\bfA_0$ may be positive constant. 
\par
The case when $\bfA_0$ is an upper or lower triangular matrix is not
covered by  Proposition~\ref{prop:2}. For example, assume that 
$\alpha_{21}=\beta_{21}=0$. Then we have the \garch\ equation
\beao
\sigma_{2,t}^2=\alpha_{02}+(\alpha_{22}Z_{2,t-1}^2+\beta_{22})\,
\sigma_{2,t-1}^2\,,\qquad t\in\bbz\,,
\eeao
which can be solved and, under the conditions of 
Proposition~\ref{prop:1}, has tail index $\alpha_2/2>0$. Writing 
$C_t=\a_{01}+(\alpha_{12} Z_{2,t-1}^2+\beta_{12})\sigma_{2,t-1}^2$, we get
\beao
\sigma_{1,t}^2= C_t+(\alpha_{11}
Z_{1,t-1}^2+\beta_{11})\,\sigma_{1,t-1}^2 \,,\qquad t\in\bbz\,.
\eeao
This is again a 1-dimensional recurrence equation but now the coefficients\\
$(C_t,\alpha_{11}
Z_{1,t-1}^2+\beta_{11})$, $t\in\bbz$, constitute a dependent strictly stationary
\seq . 
Appealing to Brandt \cite{brandt:1986}, a unique causal 
solution to this equation exists but its theoretical properties 
are not straightforward due to the 
dependence of the coefficient \seq . However, the tail of $\sigma_{1,0}^2$ is
\asy ally at least as heavy as the tail of $\sigma_{2,0}^2$. Indeed,
as $\xto$,
\beao
\P(\sigma_{1,t}^2>x)&\ge& \P(C_t>x)\\&\ge& \P((\alpha_{12}
Z_{2,t-1}^2+\beta_{12})\sigma_{2,t-1}^2>x)\\ & \sim & \E\big[(\alpha_{12}
Z_{2,t-1}^2+\beta_{12})^{\alpha_2/2}\big]\,\P(\sigma_{2,t}^2>x)\,.
\eeao
In the last step we applied  Breiman's theorem and used stationarity. 
\ere
\section{The extremogram and cross-extremogram for a bivariate
  \garch\ process}
Davis and Mikosch \cite{davis:mikosch:2009} showed for a univariate \garch\
process under the conditions of Proposition~\ref{prop:1} that
\beao
\rho_\sigma(h)&=& \lim_{\xto}\P(\sigma_h>x\mid \sigma_0>x)\\&=&
\lim_{\xto}\P(\sigma_h^2>x\mid \sigma_0^2>x)
=\E [\min(1, \Pi_h^{\alpha/2})]\,,\qquad h\ge 1\,,
\eeao 
where $A_t= \alpha_1Z_{t-1}^2+\beta_1$, $t\in\bbz$.
While the value of these quantities is not known it is possible to
determine their \asy\ order for large $h$. 
By convexity of the \fct\ $g(h)=\E[ A_0^h]$ and since $g(\alpha/2)=1$ we have
$\E [A_0^p]<1$ for
$p<\alpha/2$. Therefore
\beao
\rho_\sigma (h) \le  \E [
\min(1,  \Pi_h^p)]
\le \E [\Pi_h^p]
= (\E [A_0^p])^h\,, \qquad h\ge 1\,,
\eeao
and the \rhs\ converges to zero exponentially fast. The extremogram of
the $X$-\seq\ inherits this rate. 
Since $R_h^{0.5} 
  Z_h^+$ has a finite $\alpha$th moment, 
multiple use of Breiman's result yields
\beao
\rho_X(h)&=& 
\lim_{\xto}\P(\sigma_hZ_h>x\mid \sigma_0Z_0>x)\\
&=& \lim_{\xto}\dfrac{\P(\min(\sigma_h
  Z_h^+,\sigma_0 \,Z_0^+)>x)}{\P(\sigma_0 Z_0^+>x)}\\
&\le  & \limsup_{\xto}\dfrac{\P(\sigma_0\,\min(\Pi_h^{0.5} 
  Z_h^+,Z_0^+)> x/2)}{\P(\sigma_0 Z_0^+>x)}
+\limsup_{\xto}\dfrac{\P(R_h^{0.5} 
  Z_h^+> x/2)}{\P(\sigma_0 Z_0^+>x)}\\
&=&{\rm const}\,\dfrac{\E\big[\big( 
\min( \Pi_h^{0.5} Z_h^+,Z_0^+)\big)^{\alpha}\big]}{\E[(Z_0^+)^\alpha]}\\
&\le &{\rm const}\, (\E [A_0^p])^h\,.
\eeao
\par
Similar calculations can be done in the bivariate case.
 We restrict ourselves to the $\sigma$-\seq s.
 We assume 
the conditions of Proposition~\ref{prop:2}; in this case both
components $\sigma_{i,t}^2$, $i=1,2$, of the vector $\bfW_t$ 
in \eqref{eq:jan6b} have the same tail index.
Using relation \eqref{eq:20}, we see that
\beao
\rho_{ij}(h)&=& \lim_{\xto} \P(\sigma_{j,h}^2>x \mid
\sigma_{i,0}^2>x)\\
&=&\lim_{\xto} \dfrac{\P(\sigma_{j,h}^2>x ,
\sigma_{i,0}^2>x)}{\P(\sigma_{i,0}^2>x)}\\
&\le &\lim_{\xto} \dfrac{\P(|\bfW_h|>x , | \bfW_0|>x )}
{\P(|\bfW_0|>x) }\times 
\dfrac{\P(|\bfW_0|>x)}{\P(\sigma_{i,0}^2>x)}\,.
\eeao 
The limit of the latter ratio converges to a constant by virtue of
\regvar . 
Thus the extremograms $\rho_{ij}$ are bounded by the extremogram 
$\rho_{|\bfW|}$ of 
$(|\bfW_t|)$ times  this constant. However, \eqref{eq:20} and the
independence of $\bfW_0$ and $\bfR_h$ imply 
that for $p<\alpha/2$,
\beao
\rho_{|\bfW|}(h)&=& \lim_{\xto} \dfrac{\P(|\bfW_h|>x , | \bfW_0|>x )}
{\P(|\bfW_0|>x) }\\
&\le &\limsup_{\xto} \dfrac{\P(\|{\mathbf \Pi}_h\|\, |\bfW_0|> x/2 \,, | \bfW_0|>x )}
{\P(|\bfW_0|>x) }
+ \lim_{\xto} \P(|\bfR_h|> x/2 )\\
&=& \E [\min (1, \|{\mathbf \Pi}_h\|^{\alpha/2})]\\
&\le &\E [\min (1, \|{\mathbf \Pi}_h\|^p)]\\
&\le &  \E[ \|{\mathbf \Pi}_h\|^p]\,,\qquad h\ge 1\,.
\eeao
The \rhs\ converges to zero at an exponential rate in view of $\E
[\|{\mathbf \Pi}_{h_0}\|^p]<1$ for a sufficiently large $h_0$. 

 \section{An empirical study of the extremogram and the
  cross-extremogram}
\subsection{Estimation of the extremogram and cross-extremogram}
Davis and Mikosch \cite{davis:mikosch:2009} and Davis et al.
\cite{davis:mikosch:cribben:2012} proposed estimators of the
quantities $\rho_{ij}(h)\,,h\in\bbz\,,$ for given sets $A,B$ bounded away from zero:
\beam \label{def:empirho}
 \wh \rho _{ij}(h) = 
\dfrac{\sum_{t=1}^{n-h} \1 (X_{j,t+h}\in
 F_{X_{j,0}}^\leftarrow (1-1/m)\times B, X_{i,t} \in F_{X_{i,0}}^\leftarrow
 (1-1/m)\times   A)}{ \sum_{t=1}^n \1( X_{i,t} \in F_{X_{i,0}}^\leftarrow
 (1-1/m)   \times A)}
\eeam
for some \seq\ $m=m_n\to  \infty$ \st\ $m=o(n)$ as $\nto$.
In order to ensure standard \asy\ properties such as consistency and \asy\
normality, \cite{davis:mikosch:2009,davis:mikosch:cribben:2012}
assumed the strong mixing condition and \regvar\ for the \seq\ $(\bfX_t)$, possibly after a monotone transformation of
its components as explained in Section~\ref{subsec:regvar}. The
aforementioned growth conditions on the \seq\ $(m_n)$ are standard in
extreme value statistics and cannot be avoided. They ensure that
sufficiently high thresholds $F_{X_{i,0}}^\leftarrow (1-1/m)$,
$i=1,2,$ are chosen. These thresholds guarantee that a certain
fraction of the data is taken which may be considered extreme as
regards their distance from the origin. For practical purposes, we 
take the corresponding empirical $(1-1/m)$-quantiles of the
components.
\par
Although central limit theory can be shown for $\wh \rho _{ij}$ at a
finite number of lags $h$, the \asy\ covariance structure is not
tractable. Davis et al. \cite{davis:mikosch:cribben:2012} propose two
methods for the construction of credible confidence bands: the
stationary bootstrap and random permutations. In this paper, we stick
to the latter procedure. It is based on the simple idea that, if the sample
$\bfX_1,\ldots,\bfX_n$ were iid, 
random permutations of the sample would not change its dependence
structure, hence the extremogram and cross-extremogram would not
change. In what follows, we calculate the (cross-) extremograms based on 100 
random permutations of the sample. First we calculate the 100
extremogram values at each lag. Then we choose the 96\% empirical
quantile at each lag and finally take the maximum over the lags of
interest. This value is shown as a solid horizontal line in the corresponding graphs. 
This procedure is quick and clean: if the sample (cross-) extremogram
at a given lag is above the horizontal line
this is an indication of disagreement with the iid hypothesis.

\subsection{Simulated \garch\ data}
We provide a brief study of the sample (cross)-extremograms 
of simulated bivariate \garch\ processes and their residuals. 
We choose bivariate \garch\ models with iid bivariate $t$-distributed innovations  
$(\bfZ_t)$ with $10$ degrees of freedom and covariance matrix $P$ given in \eqref{eq:14}.  
We simulate from the model \eqref{eq:15} with parameters $(\alpha_{01},\alpha_{02})=(10^{-6},10^{-6})$
(the magnitude of these parameters is in agreement with values estimated from return data) and specified matrices and correlations 
\begin{align}
\label{para-alpha-beta}
\left(
    \begin{array}{cc}
      \alpha_{11} & \alpha_{12}  \\
      \alpha_{12} & \alpha_{22} 
    \end{array}
  \right),\quad \left(
    \begin{array}{cc}
      \beta_{11} & \beta_{12}  \\
      \beta_{12} & \beta_{22} 
    \end{array}
  \right)\,,\quad \rho\,.
\end{align}
We start by considering examples with $\rho=0$ (Examples (1)--(4))
and $\rho=0.7$ (Examples (5)--(6)) with respective symmetric parameter
matrices \eqref{para-alpha-beta}:
\begin{align*}
&(1)\quad
\left(
    \begin{array}{cc}
      .1 & 0  \\
       0 & .1 
    \end{array}
  \right), \left(
    \begin{array}{cc}
      .8 & 0  \\
       0 & .8
    \end{array}
  \right),\quad 
(2)\, 
\left(
    \begin{array}{cc}
      .1 & .05  \\
       .05 & .1 
    \end{array}
  \right), \left(
    \begin{array}{cc}
      .8 & 0  \\
       0 & .8
    \end{array}
  \right),\quad
(3)\, 
\left(
    \begin{array}{cc}
      .1 & 0  \\
       0 & .1 
    \end{array}
  \right), \left(
    \begin{array}{cc}
      .8 & .04  \\
       .04 & .8
    \end{array}
  \right), \\
(4)&\quad
\left(
    \begin{array}{cc}
      .1 & .02  \\
       .02 & .1 
    \end{array}
  \right), \left(
    \begin{array}{cc}
      .8 & .04  \\
       .04 & .8
    \end{array}
  \right),\quad 
(5)\, 
\left(
    \begin{array}{cc}
      .1 & 0  \\
       0 & .1 
    \end{array}
  \right), \left(
    \begin{array}{cc}
      .8 & 0  \\
       0 & .8
    \end{array}
  \right),\quad
(6)\, 
\left(
    \begin{array}{cc}
      .1 & .02  \\
      .02 & .1 
    \end{array}
  \right), \left(
    \begin{array}{cc}
      .8 & .04  \\
       .04 & .8
    \end{array}
  \right). 
\end{align*}
Here we always choose small $\alpha$-values while the diagonal $\beta$-values are close to 1. This is in agreement
with estimated parameters on return data.  
We generate samples of size $n=50,000$, using the R package 'ccgarch'\footnote{
Note that estimation with  ``ccgarch'' requires choosing initial
values. In most cases, we first examine componentwise univariate
\garch\ fits by the R package ``fGarch'' and then we choose  these estimates as
initial values. If the univariate estimation does not converge we try several initial values on a grid of size
$0.1$. In this case, the estimates sometimes differ by attaining local
minima. Judging from the residuals, the eigenvalues of the estimated parameters 
\eqref{eq:spectral} and the values of the likelihood functions, we choose
an ``optimal'' estimator. Except for one case of stock
return data (see Section \ref{sbsec:stock}), this procedure works.}, and
calculate the (cross)-extremograms $\wh \rho_{ij}(h)$ in
\eqref{def:empirho} with $A=B=(1,\infty)$.  
\begin{figure} 
\begin{center}
\includegraphics[width=0.475\textwidth]{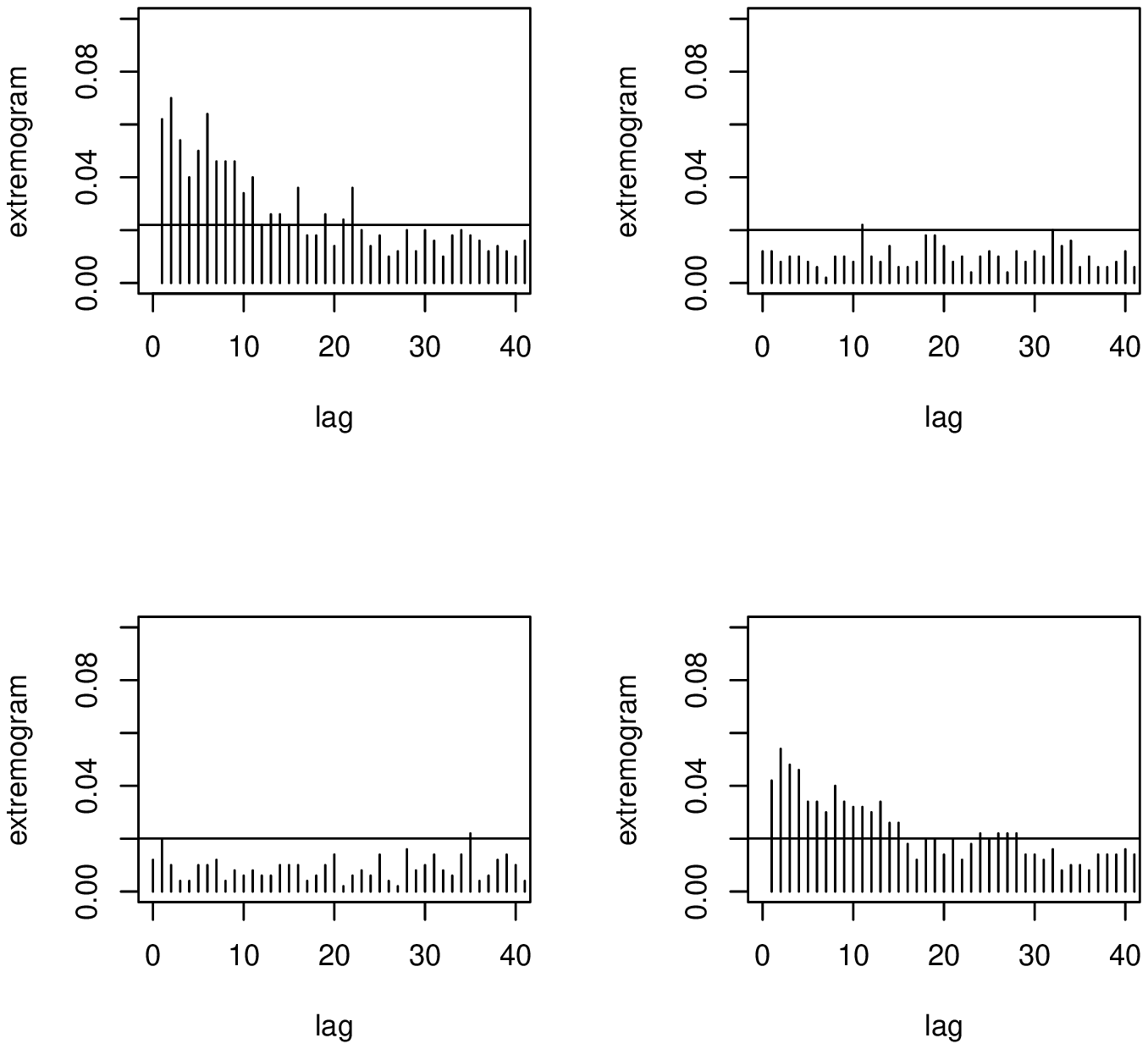}
\includegraphics[width=0.475\textwidth]{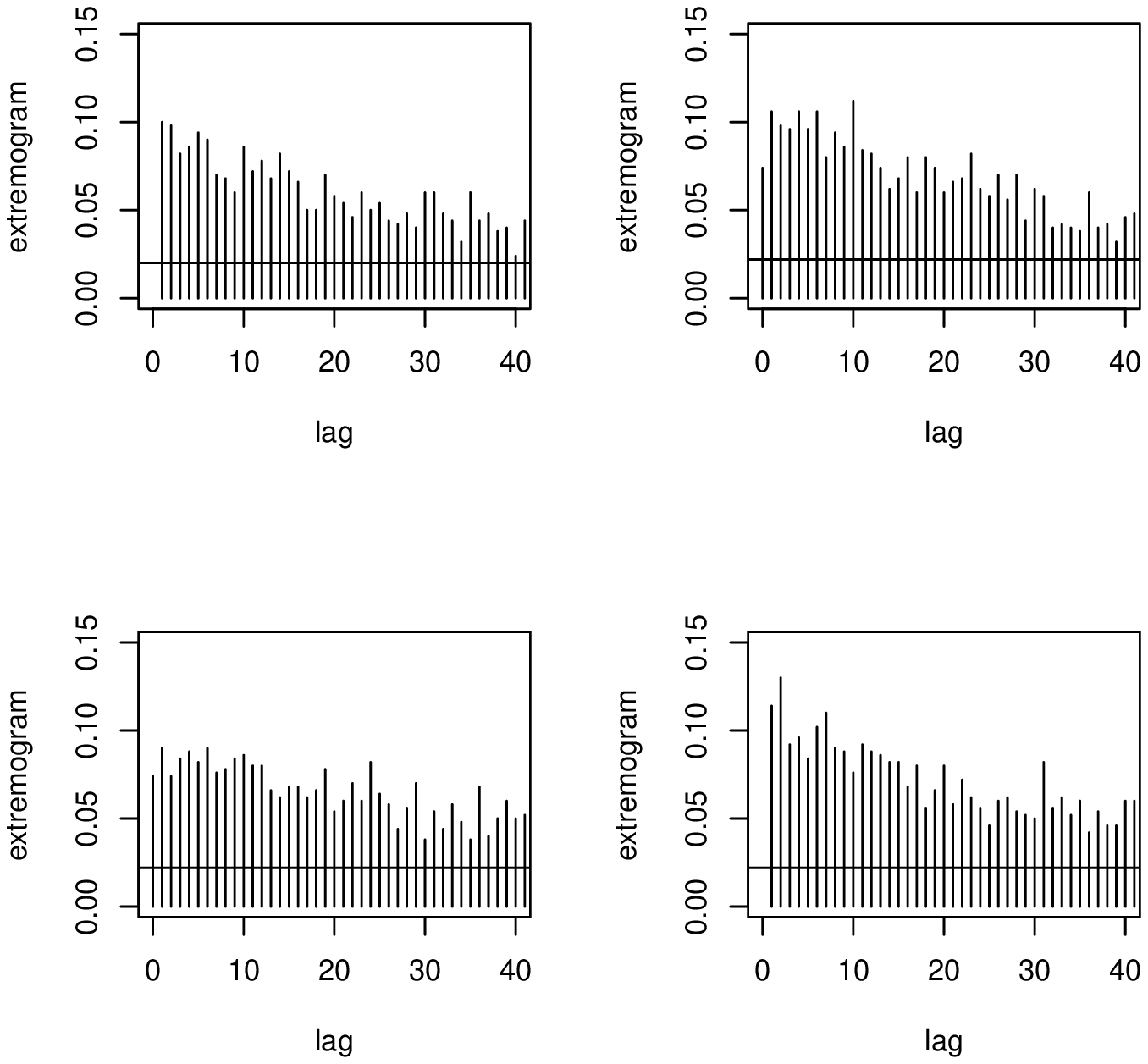}
\includegraphics[width=0.475\textwidth]{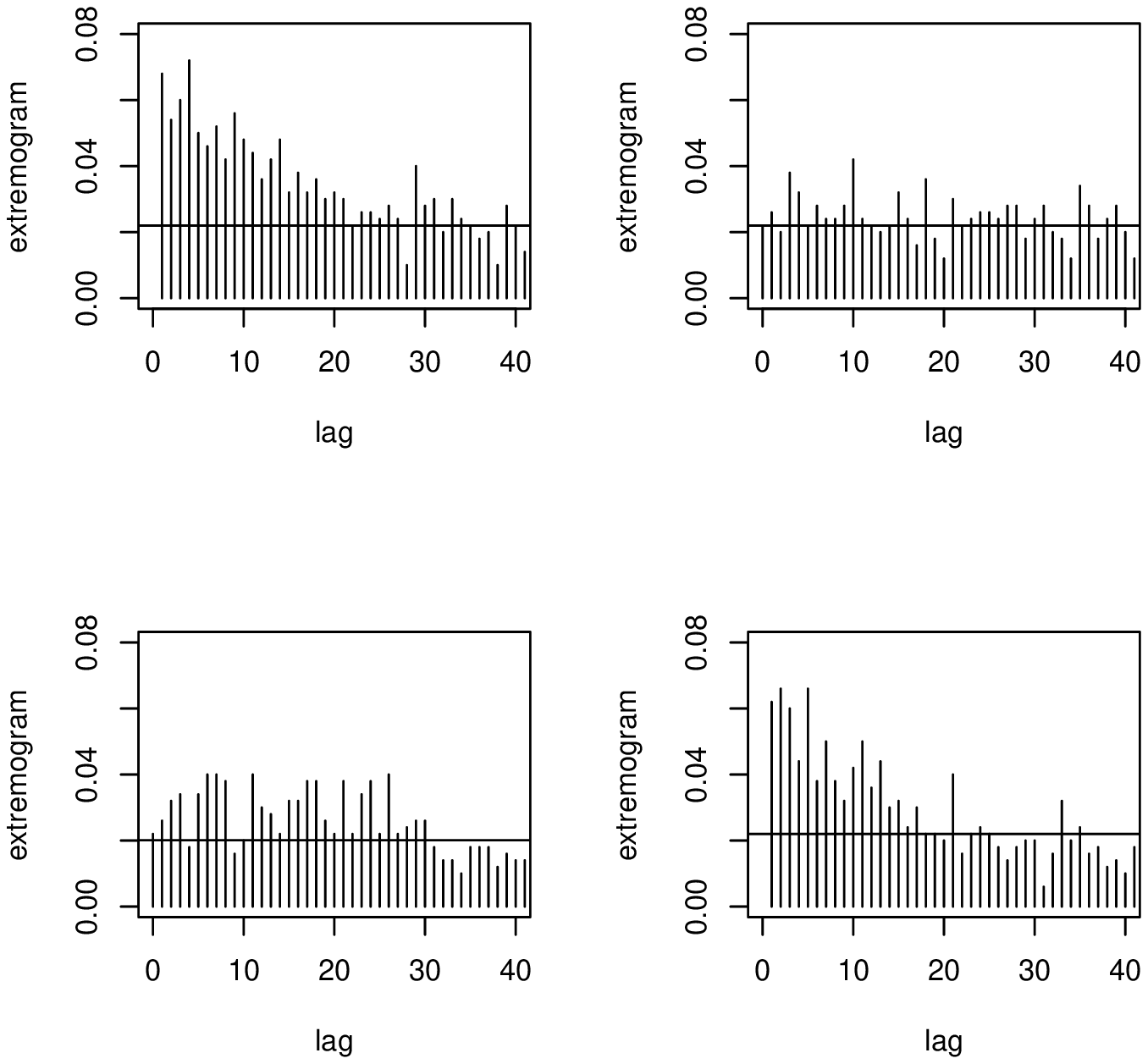}
\includegraphics[width=0.475\textwidth]{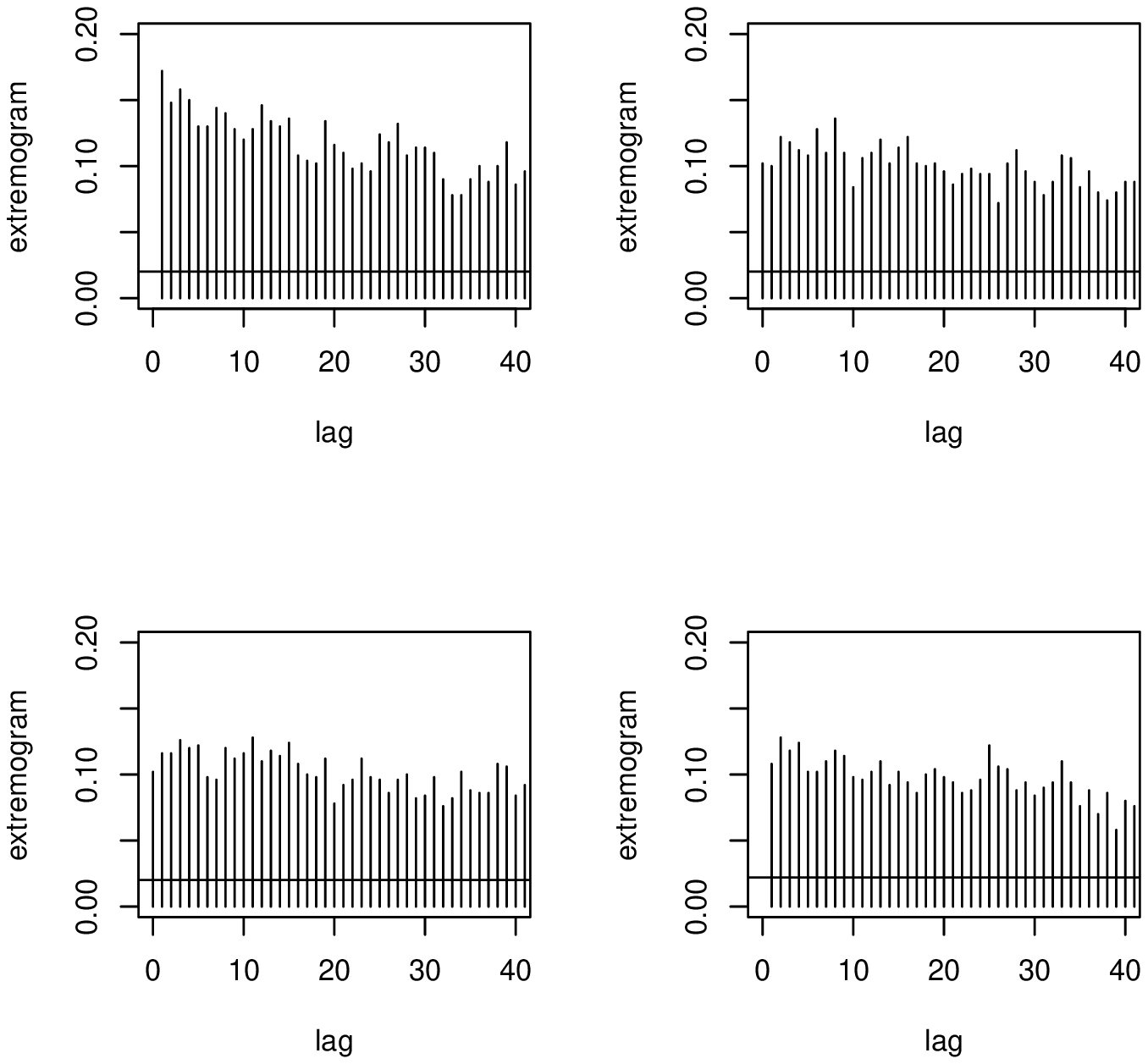}
\includegraphics[width=0.475\textwidth]{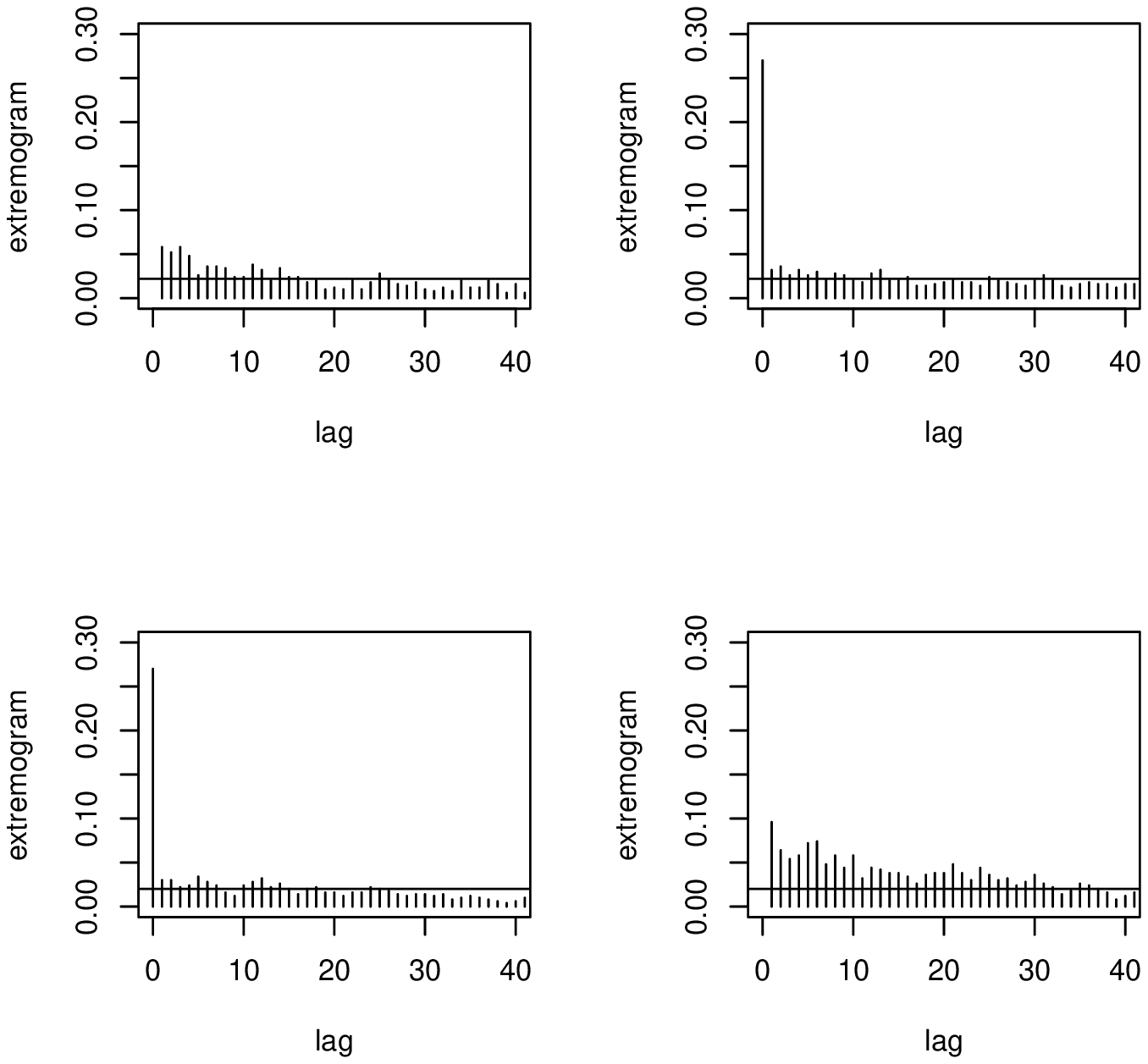}
\includegraphics[width=0.475\textwidth]{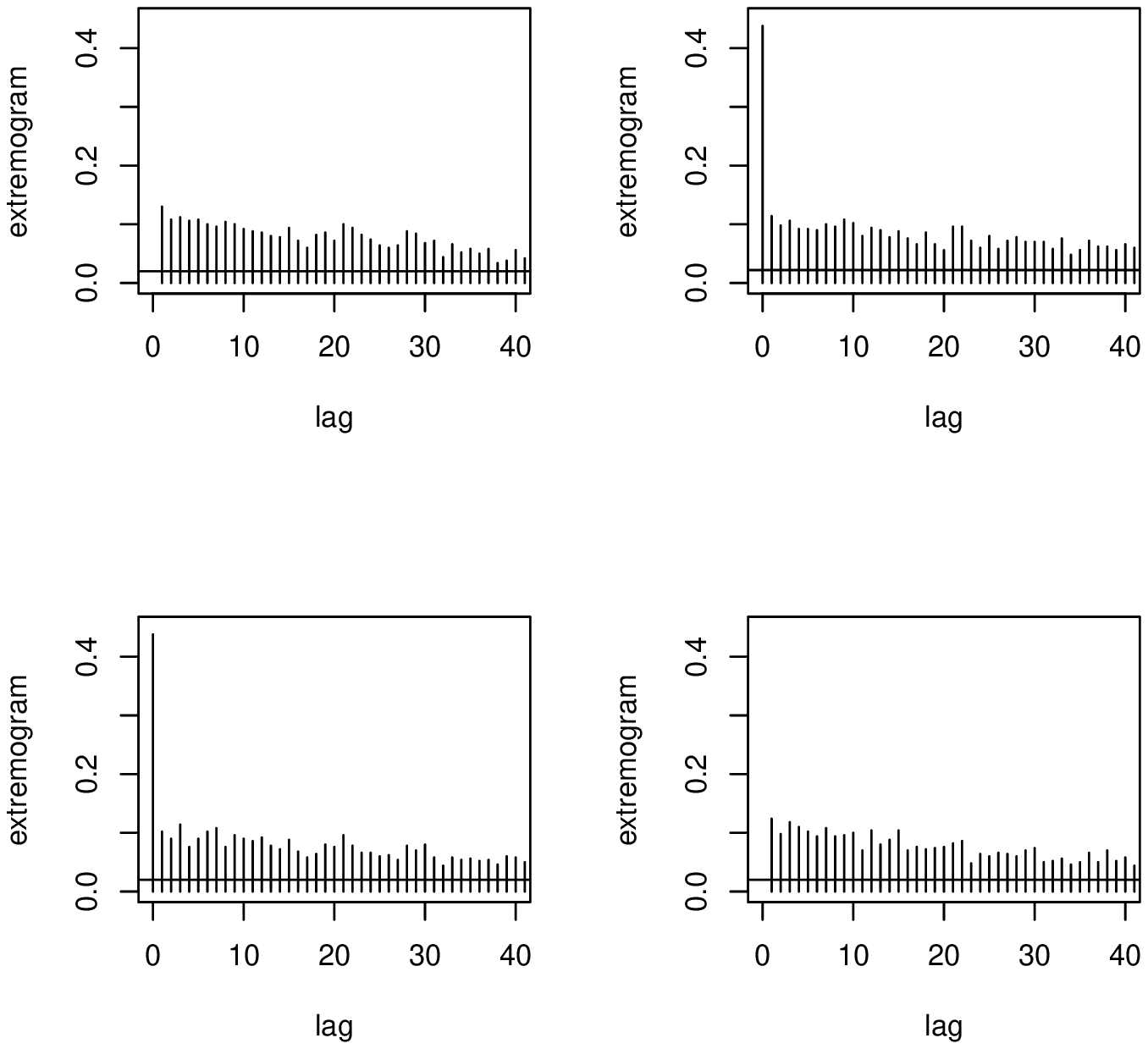}
\end{center}
\caption{(Cross-) extremograms corresponding to Examples (1) 
(top left $2\times 2$ graphs), (2) (top right), (3) (middle left), (4) (middle right), (5)
 (bottom left) and (6) (bottom right).
}
\label{fig:1-1}
\end{figure}
The simulation results for Examples (1)--(6) are given in 
Figure \ref{fig:1-1}. These figures indicate that small changes in the $\alpha$- or $\beta$-values 
lead to substantial changes in the extremal dependence structure. 
In all cases we observe serial extremal dependence in extremograms
(diagonal parts of $2 \times 2$ graphs).
In Examples (5)--(6) we also observe large spikes 
in the cross-extremograms at lag zero due to $\rho \neq 0$. This is in
contrast to Examples (1)--(4) with $\rho=0$. Example (1) shows clear asymmetry in
cross-extremograms compared with Example (2), namely, no-correlation can
be can be read in the Example (1), which reflects the componentwise
independence setting. Others exhibit dependencies in cross-extremograms
to greater or lesser degrees. 
\par
\begin{figure} 
\begin{center}
\includegraphics[width=0.475\textwidth]{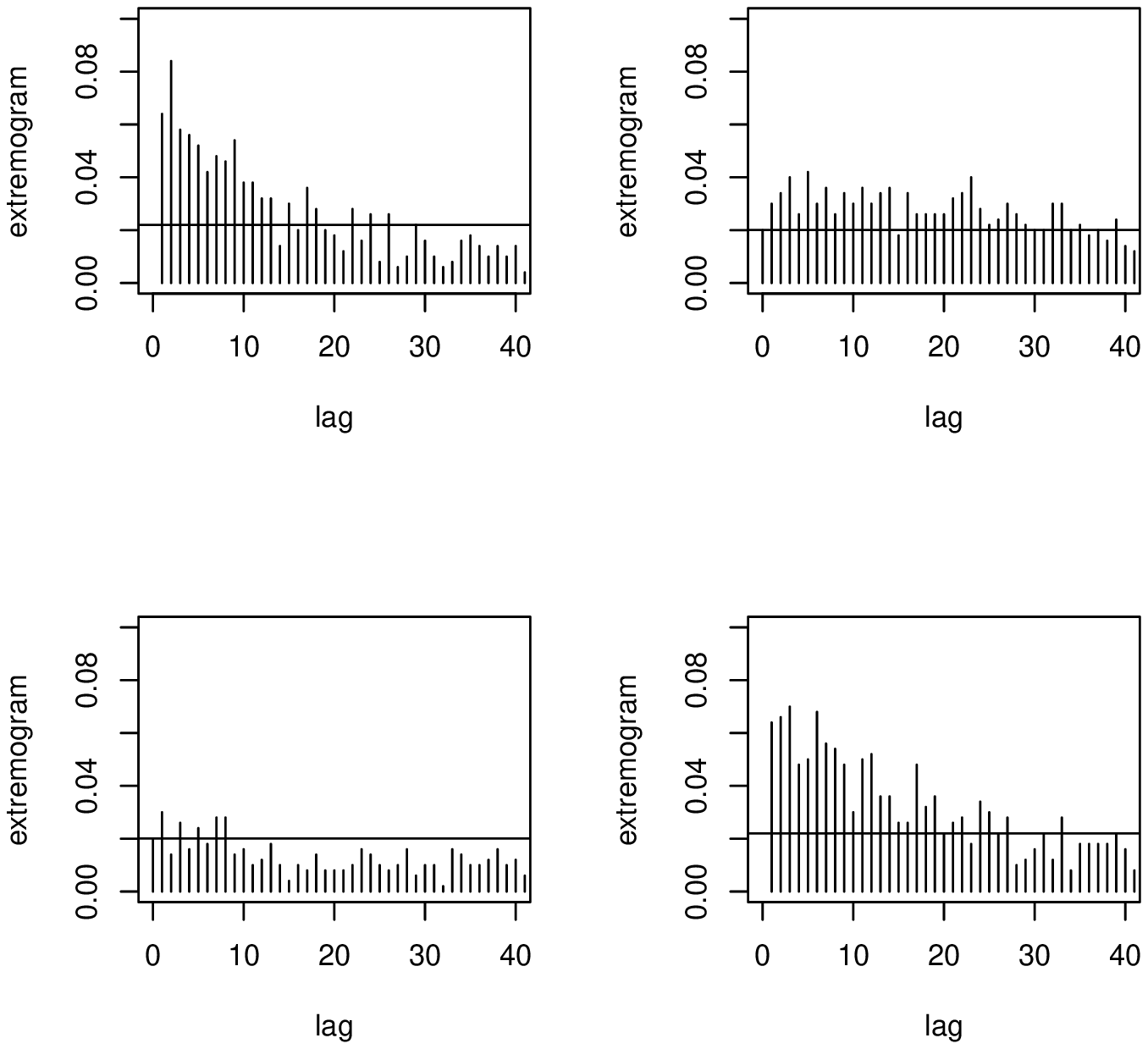}
\includegraphics[width=0.475\textwidth]{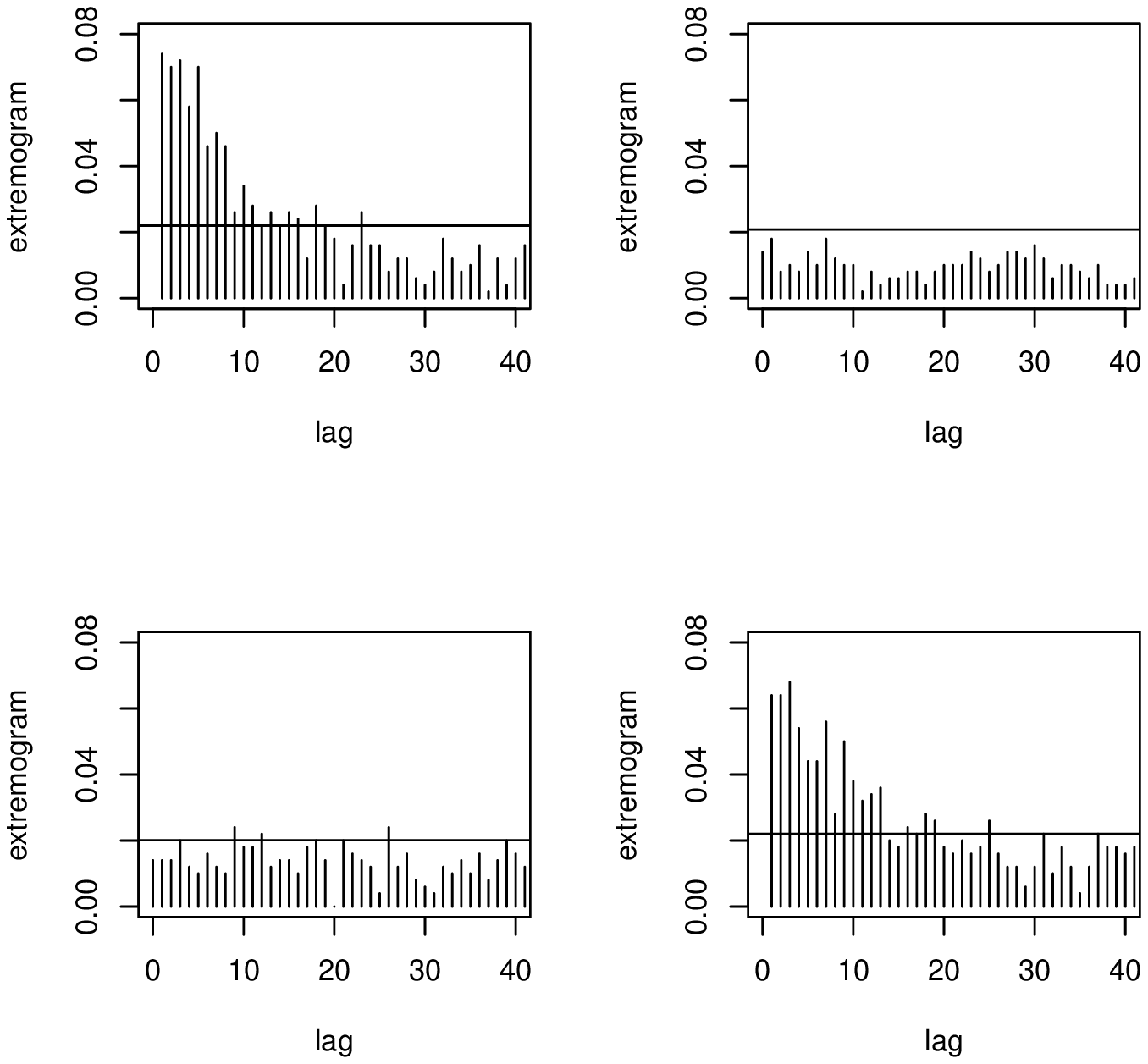}
\includegraphics[width=0.475\textwidth]{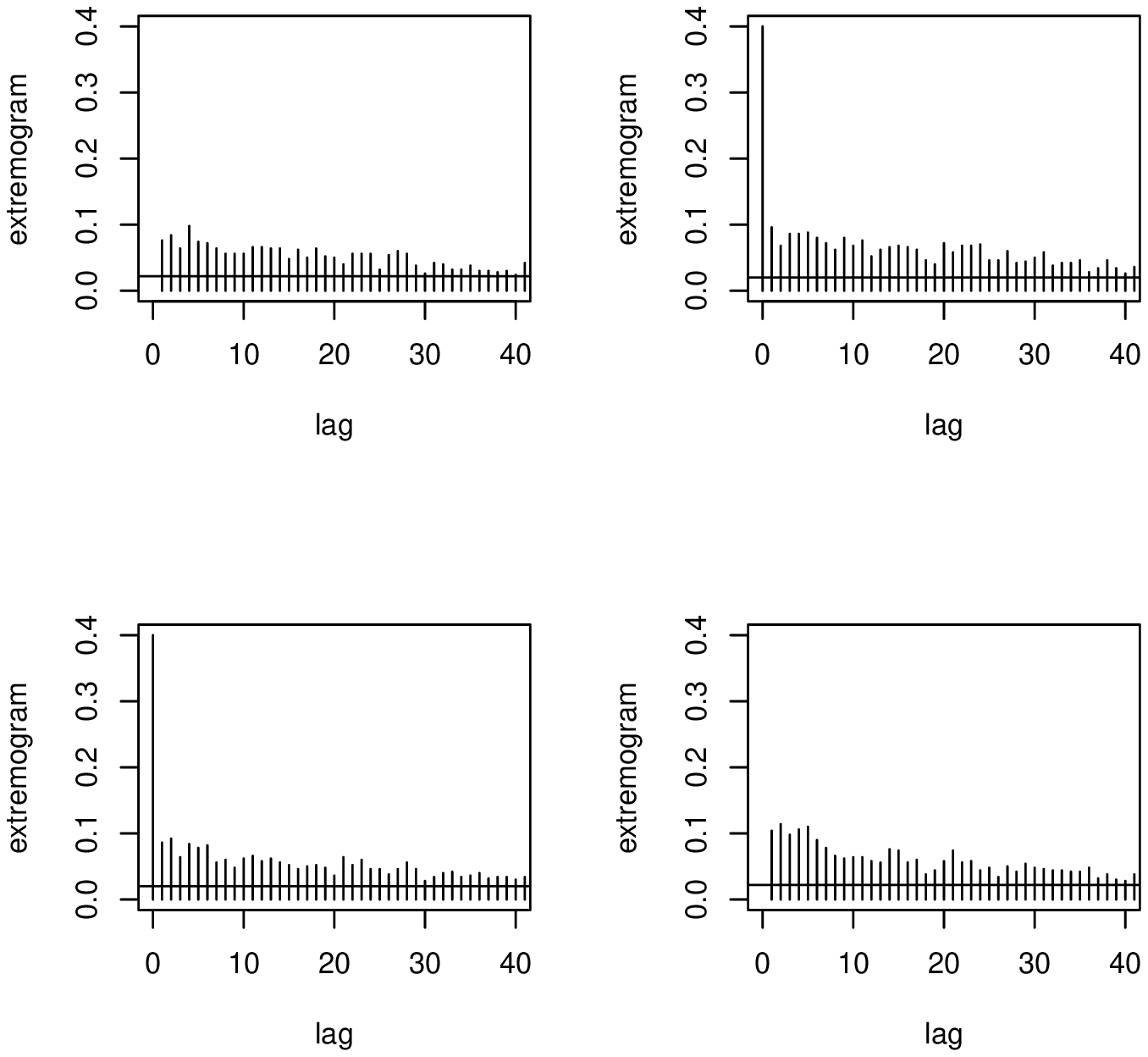}
\includegraphics[width=0.475\textwidth]{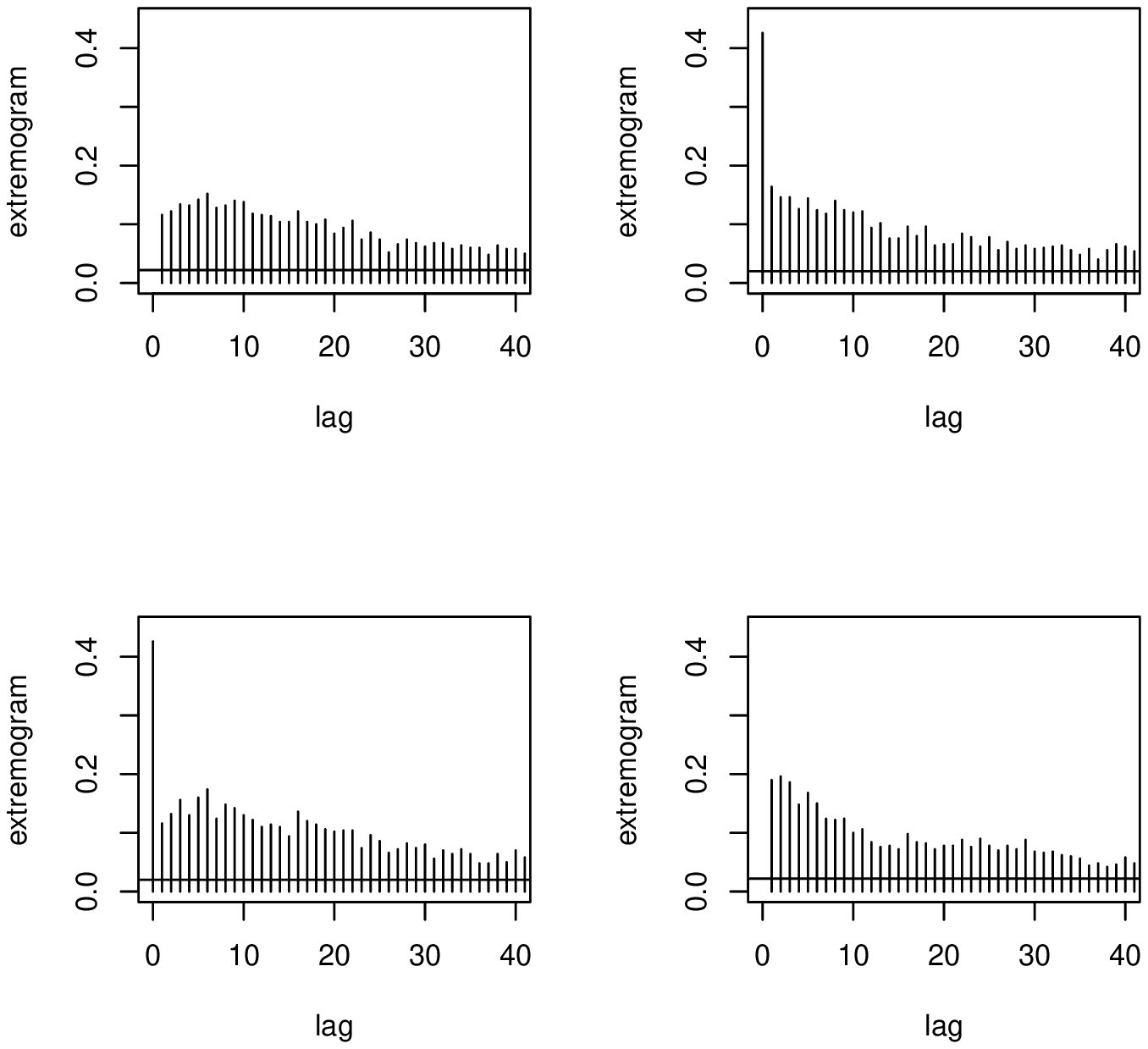}
\end{center}
\caption{(Cross-) extremograms of Examples (7) (top left $2\times 2$ graphs), 
(8) (top right), (9) (bottom left) and  (10) (bottom right).} \label{fig:1-2}
\end{figure}

In Examples (7)--(8) and (9)--(10) we choose $\rho=0$ and 
$\rho=0.7$, respectively, and the asymmetric $\alpha$- and
$\beta$-matrices \eqref{para-alpha-beta} as follows 
\begin{align*}
(7)\quad
\left(
    \begin{array}{cc}
      .1 & 0  \\
      .07 & .1 
    \end{array}
  \right), \left(
    \begin{array}{cc}
      .8 & 0  \\
       0 & .8
    \end{array}
  \right),\quad 
(8)\, 
\left(
    \begin{array}{cc}
      .1 & 0  \\
       0 & .1 
    \end{array}
  \right), \left(
    \begin{array}{cc}
      .8 & .04  \\
       0 & .8
    \end{array}
  \right), \\
(9)\quad
\left(
    \begin{array}{cc}
      .1 & 0  \\
      .05 & .1 
    \end{array}
  \right), \left(
    \begin{array}{cc}
      .8 & .04  \\
       0 & .8
    \end{array}
  \right),\quad 
(10)\, 
\left(
    \begin{array}{cc}
      .1 & 0  \\
       .1 & .2 
    \end{array}
  \right), \left(
    \begin{array}{cc}
      .8 & .07  \\
       0 & .6
    \end{array}
  \right). 
\end{align*}
Figure \ref{fig:1-2} indicates that the  
asymmetry manifests through the $\alpha$-matrix rather  than
the $\beta$-matrix. Again, we observe
large spikes at lag zero for the cross-extremograms when $\rho\neq 0$
(Examples (8) and (10)). 
Example (10) shows the effects when 
diagonal elements are distinct. 

Our next goal is to show (cross-) extremograms of the residuals of simulated bivariate \garch\ models. Although we know the 
innovations \seq\ in this case, we want to illustrate how standard MLE techniques work. Of course, we expect that
the residuals of the models have properties close to those of an iid \seq\ also as regards their extremal behavior. 
The estimation is done in two ways:  (1) we fit component-wise
univariate \garch\ models, applying  MLE and assuming
student $t$-distributions for the innovations; (2) following Ling and McAleer \cite{ling:mcaleer:2003}, we apply
bivariate Gaussian quasi-MLE (QMLE). We consider the model \eqref{eq:15} with
given parameter $(\alpha_{01},\alpha_{02})=(10^{-6},10^{-6})$ and parameter matrices 
\eqref{para-alpha-beta} as follows
\beao
\barr{lllllll}
&(11)& \left(
    \begin{array}{cc}
      .1 & 0  \\
      .05 & .1
    \end{array}
  \right),  \left(
    \begin{array}{cc}
      .8 & .03  \\
      0 & .8
    \end{array}
  \right), &\rho=.7,&
(12)&\left(
    \begin{array}{cc}
      .2 & 0  \\
      .07 & .1
    \end{array}
  \right), \left(
    \begin{array}{cc}
      .7 & 0  \\
      .02 & .5
    \end{array}
  \right),&\rho=.5\earr
\eeao
Component-wise univariate MLE yields the following results:
\begin{table}[htb]
\begin{minipage}{0.45\hsize}
\begin{center}
  \begin{tabular}{r|rrr}
   Ex. (11)   & $\hat \alpha_{i}$ & $\hat \beta_{i}$ & degree for $t$ \\ \hline
   $i=1$ & $.137$ & $.831 $ & $9.81$  \\
   $i=2$ & $.169$ & $.802$  & $10.00$
  \end{tabular}
\end{center}
\end{minipage}
\begin{minipage}{0.45\hsize}
\begin{center}
\begin{tabular}{r|rrr}
    Ex. (12) & $\hat \alpha_i$ & $\hat \beta_i$ & degree for $t$ \\ \hline
   $i=1$ & $.254$ & $.698 $ & $10.00$  \\
   $i=2$ & $.232$ & $.633$  & $7.88$ 
\end{tabular}
\end{center}
\end{minipage}
\end{table}
Despite the misspefication of a bivariate \garch\ model, the univariate estimation 
leads to reasonable estimation results except for the second component of Example (12), where the estimated parameters are far from
the true ones.
\par
Bivariate QMLE yields the following estimation results for parameters \eqref{para-alpha-beta}:
\beao
\barr{lllllll}
&(11)& \left(
    \begin{array}{cc}
      .130 & 0  \\
      .056 & .125
    \end{array}
  \right),  \left(
   \begin{array}{cc}
      .778 & .025  \\
      .039 & .790
    \end{array}
  \right), &\ .7,&
(12)&\left(
    \begin{array}{cc}
      .276 & 0  \\
      .078 & .071
    \end{array}
  \right), \left(
    \begin{array}{cc}
      .685 & .004  \\
      .004 & .638
    \end{array}
  \right),&\ .504\earr
\eeao
\begin{figure} 
\begin{center}
\includegraphics[width=0.45\textwidth]{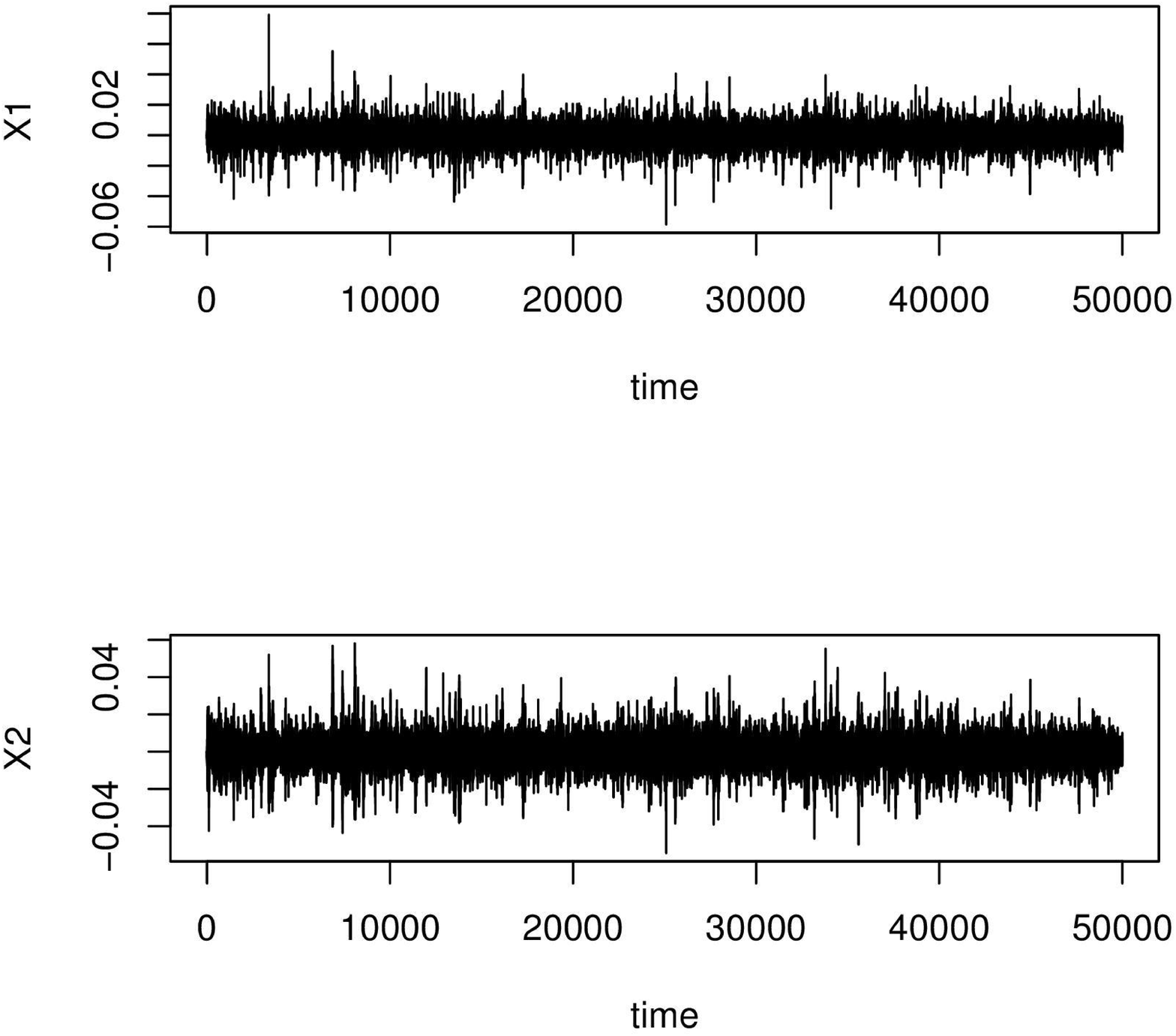}
\includegraphics[width=0.45\textwidth]{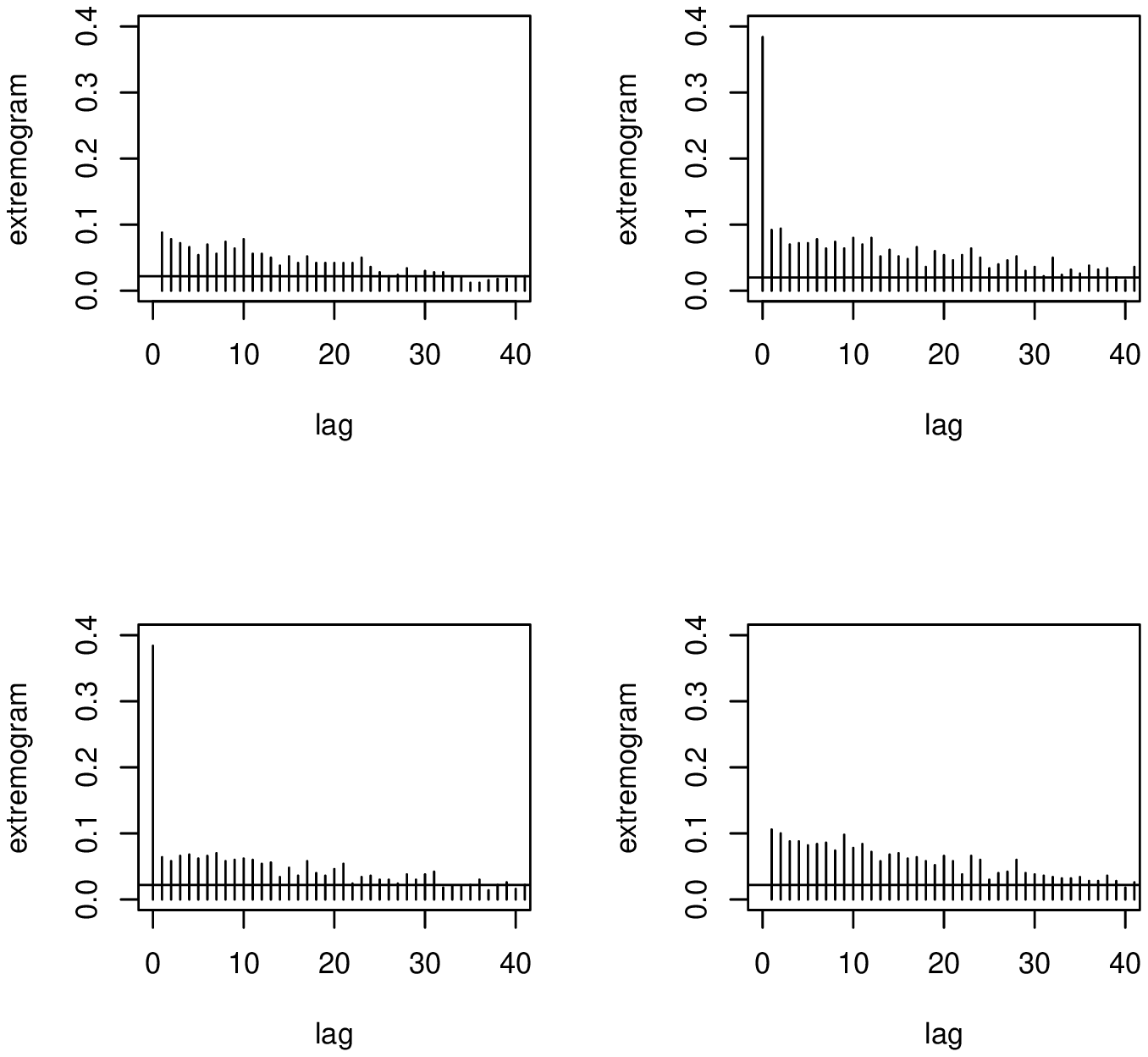} 
\includegraphics[width=0.45\textwidth]{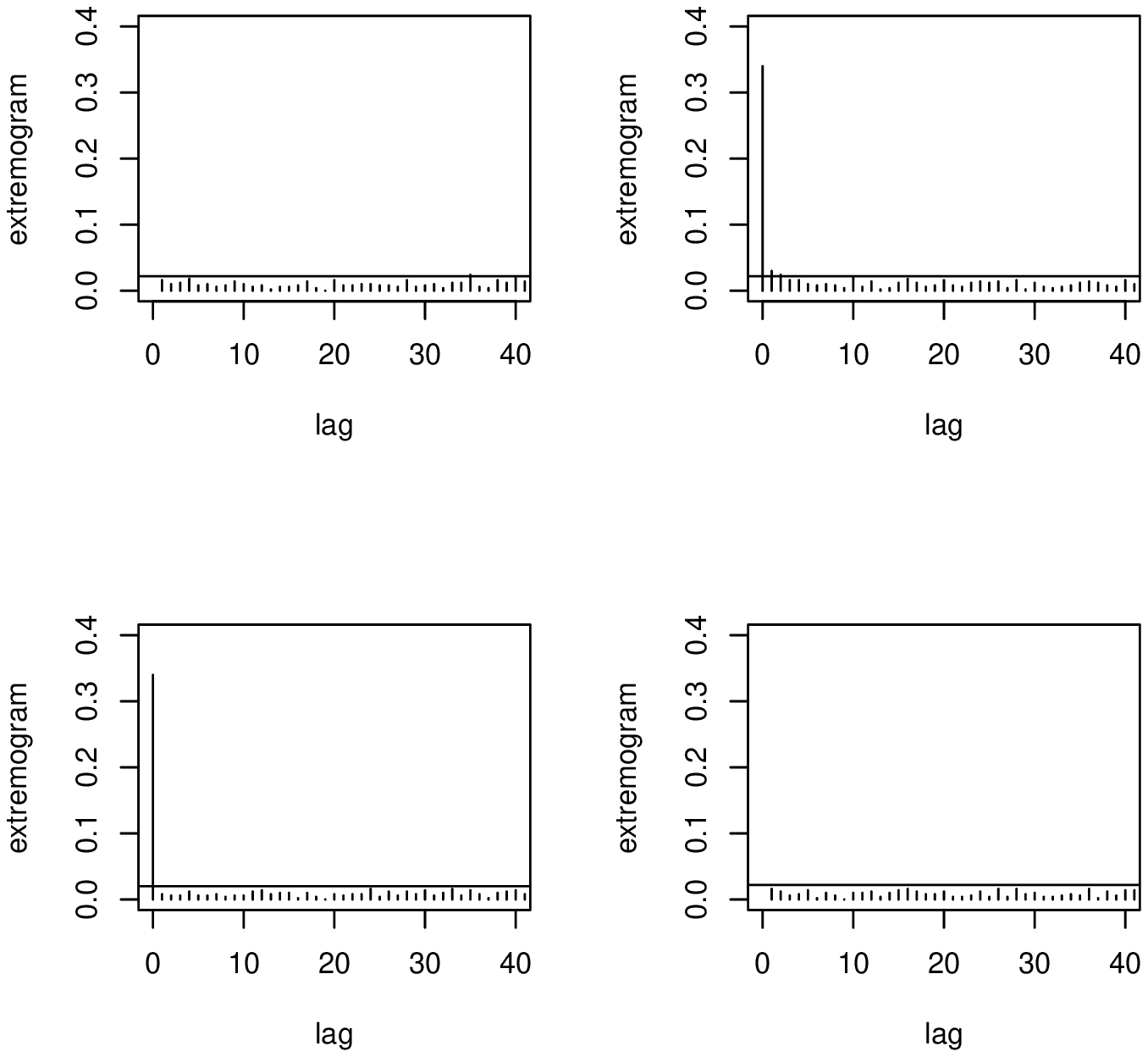}
\includegraphics[width=0.45\textwidth]{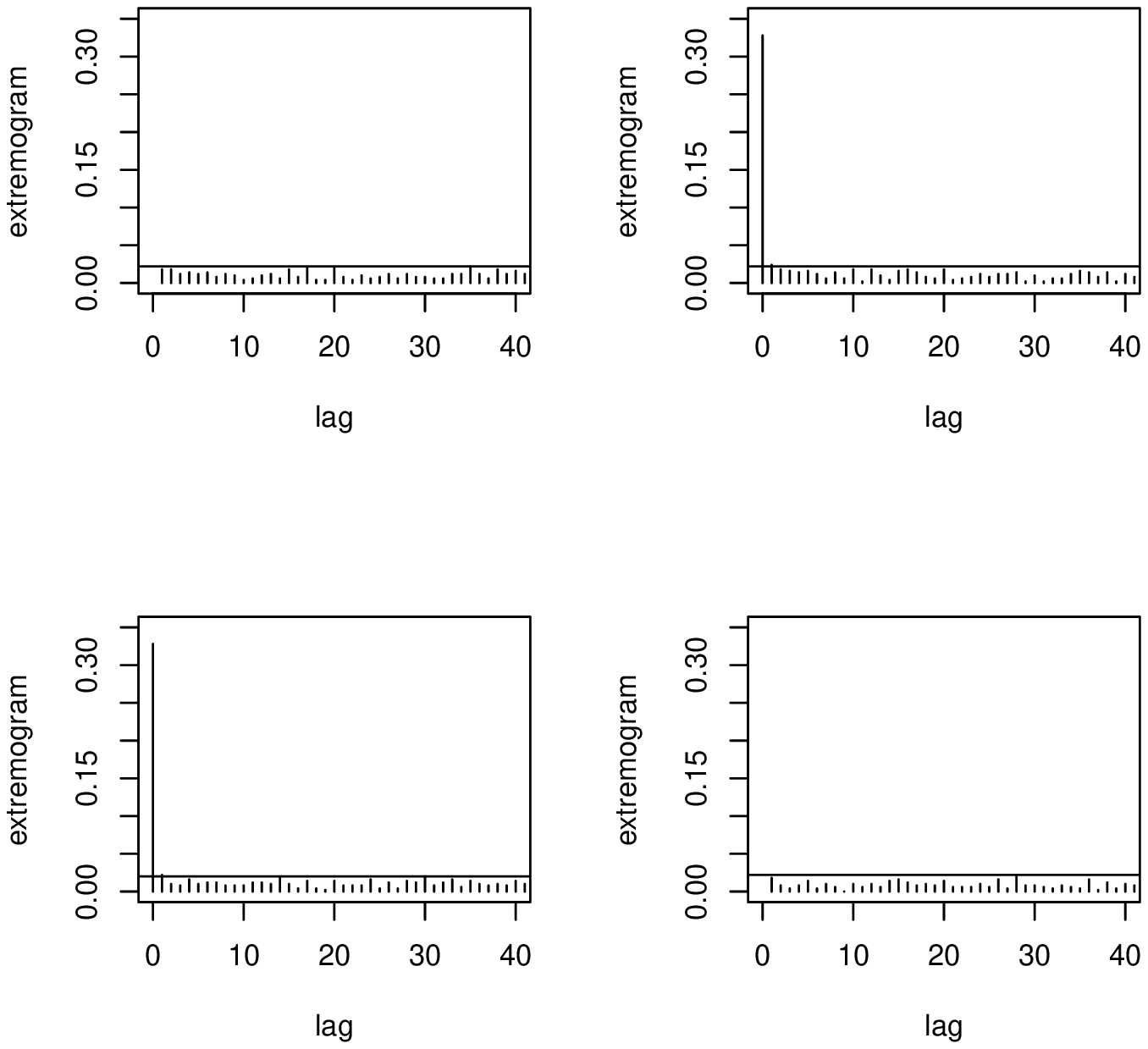}
\caption{(Cross-) extremograms in Example (11). Top left: bivariate
 \garch\ data. 
 Top right $2\times 2$ graphs: (Cross-) extremogram of bivariate \garch\ data, 
 Bottom left $2\times 2$ graphs:  
(Cross-) extremogram of residuals after a bivariate QMLE fit, 
 Bottom right $2\times 2$ graphs: (Cross-) extremogram of residuals
 after applying the componentwise MLE.}
\label{fig:1-11}
\end{center}
\end{figure}
\begin{figure} 
\begin{center}
\includegraphics[width=0.45\textwidth]{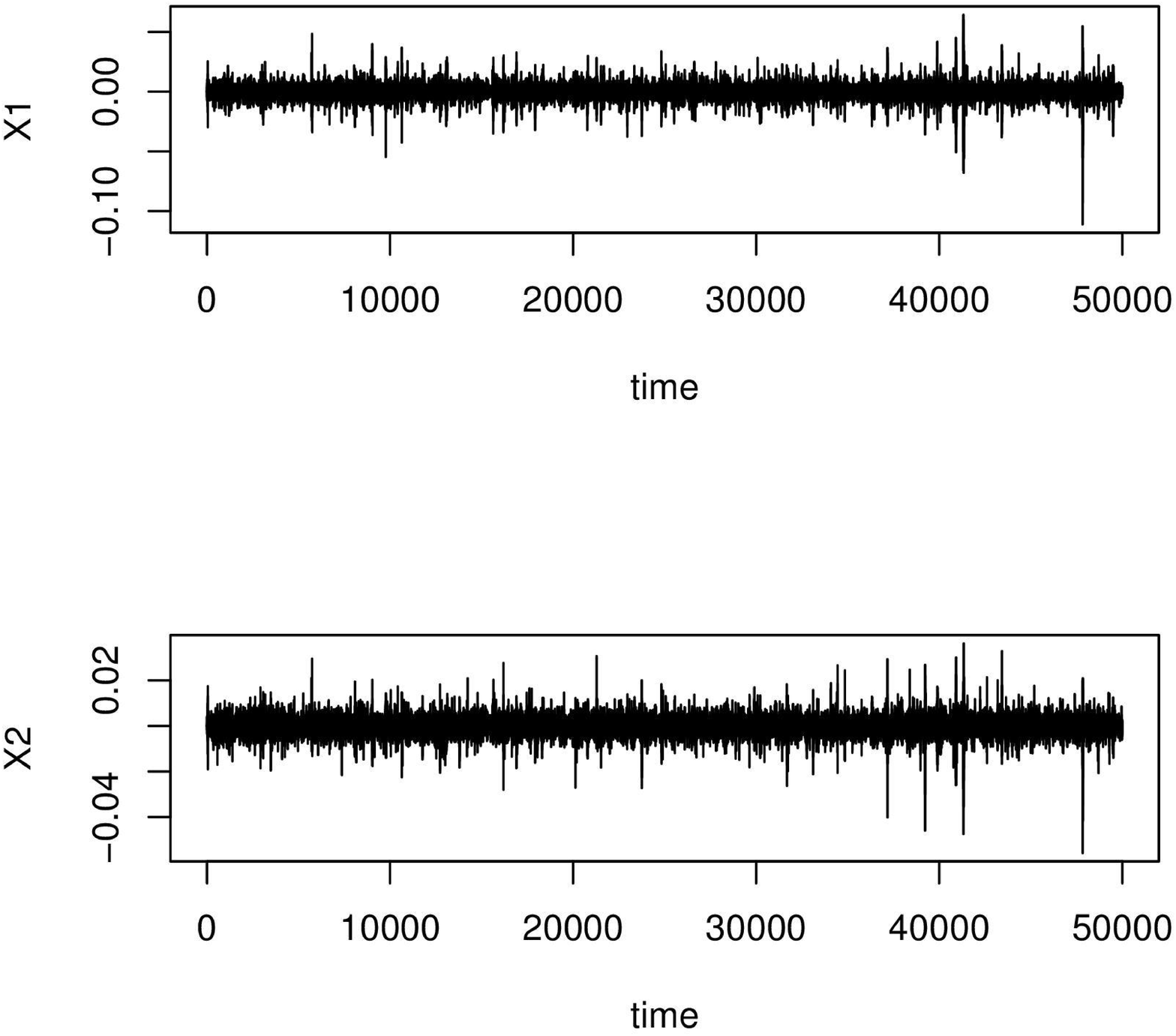}
\includegraphics[width=0.45\textwidth]{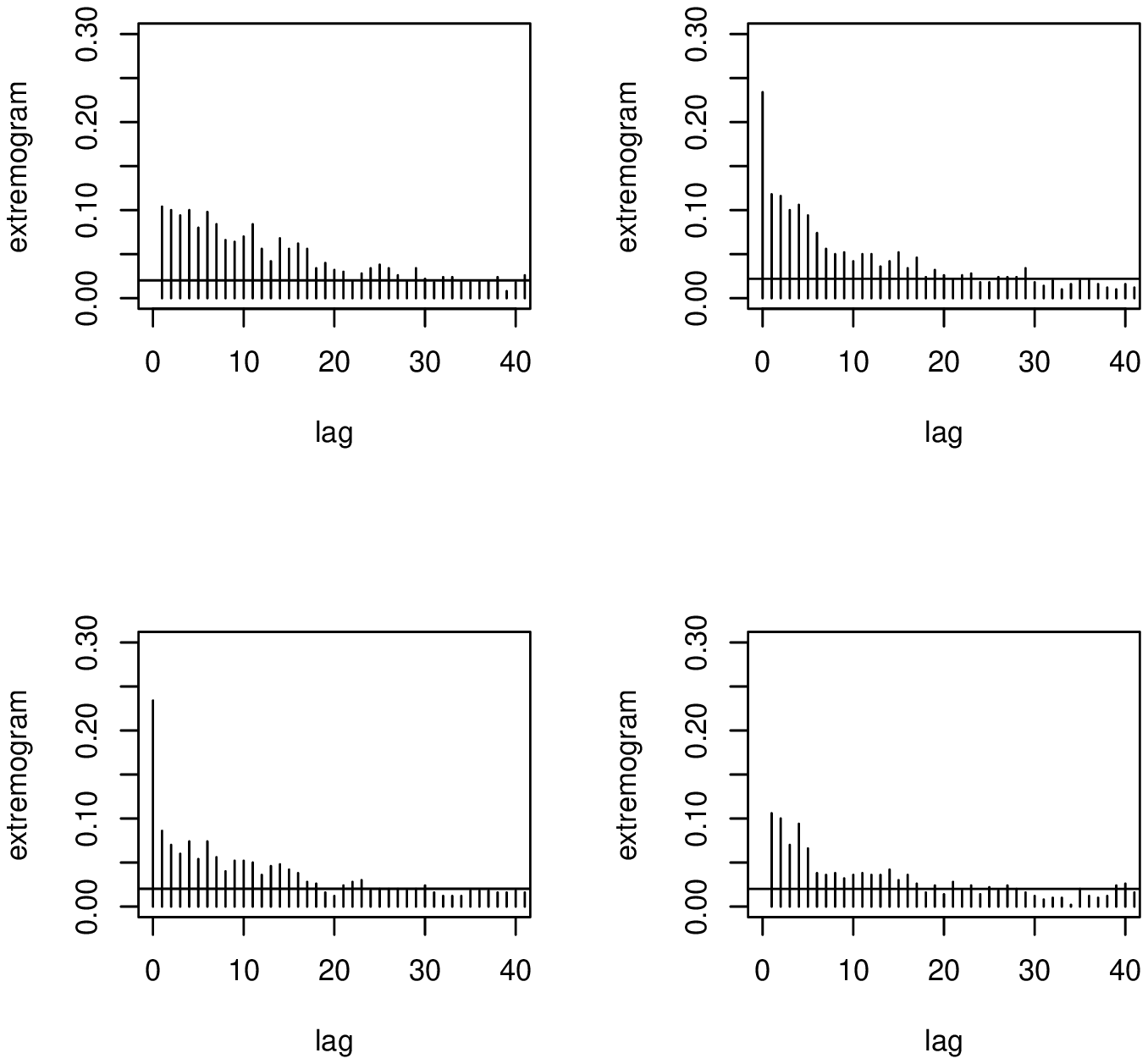}
\includegraphics[width=0.45\textwidth]{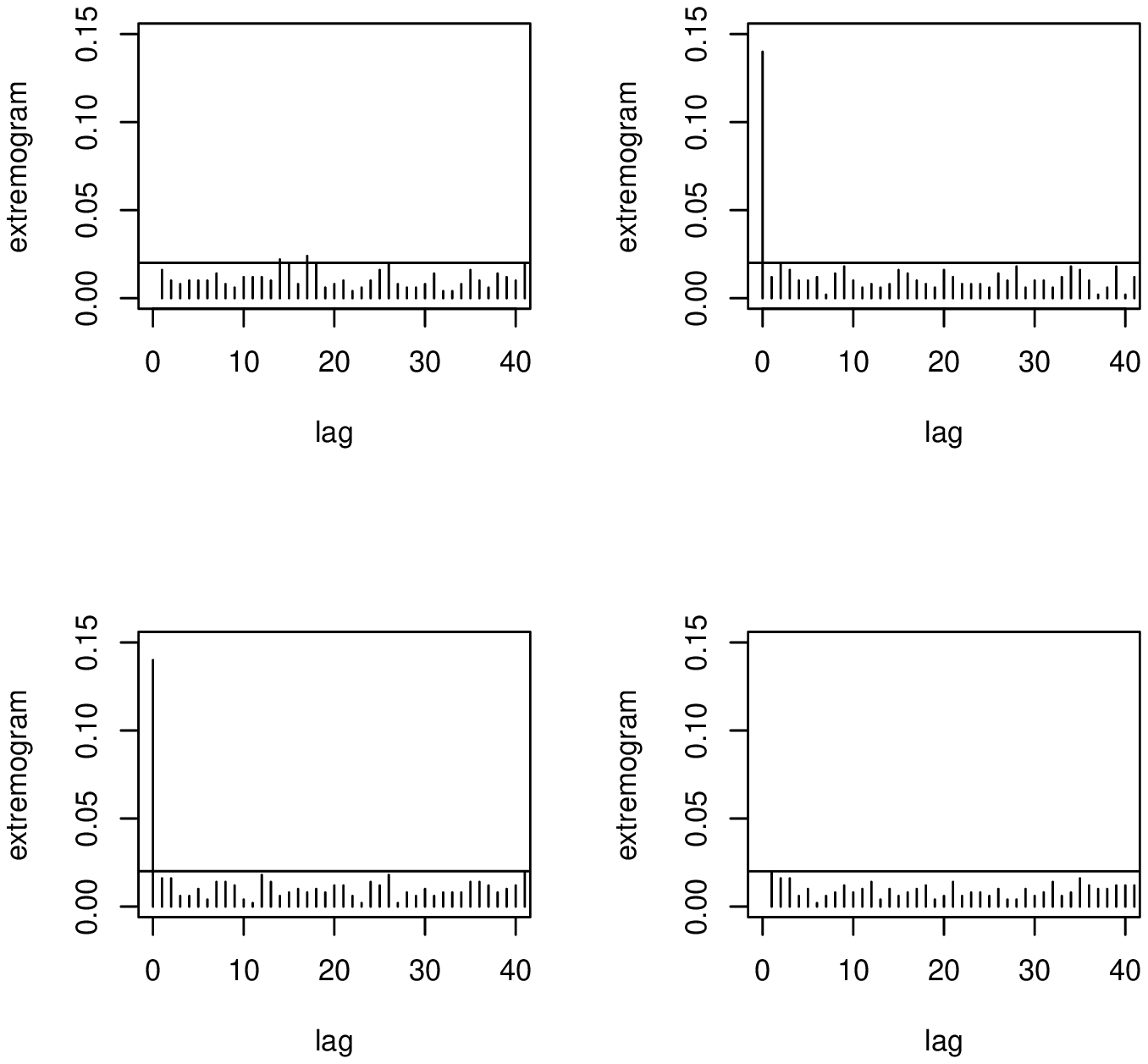}
\includegraphics[width=0.45\textwidth]{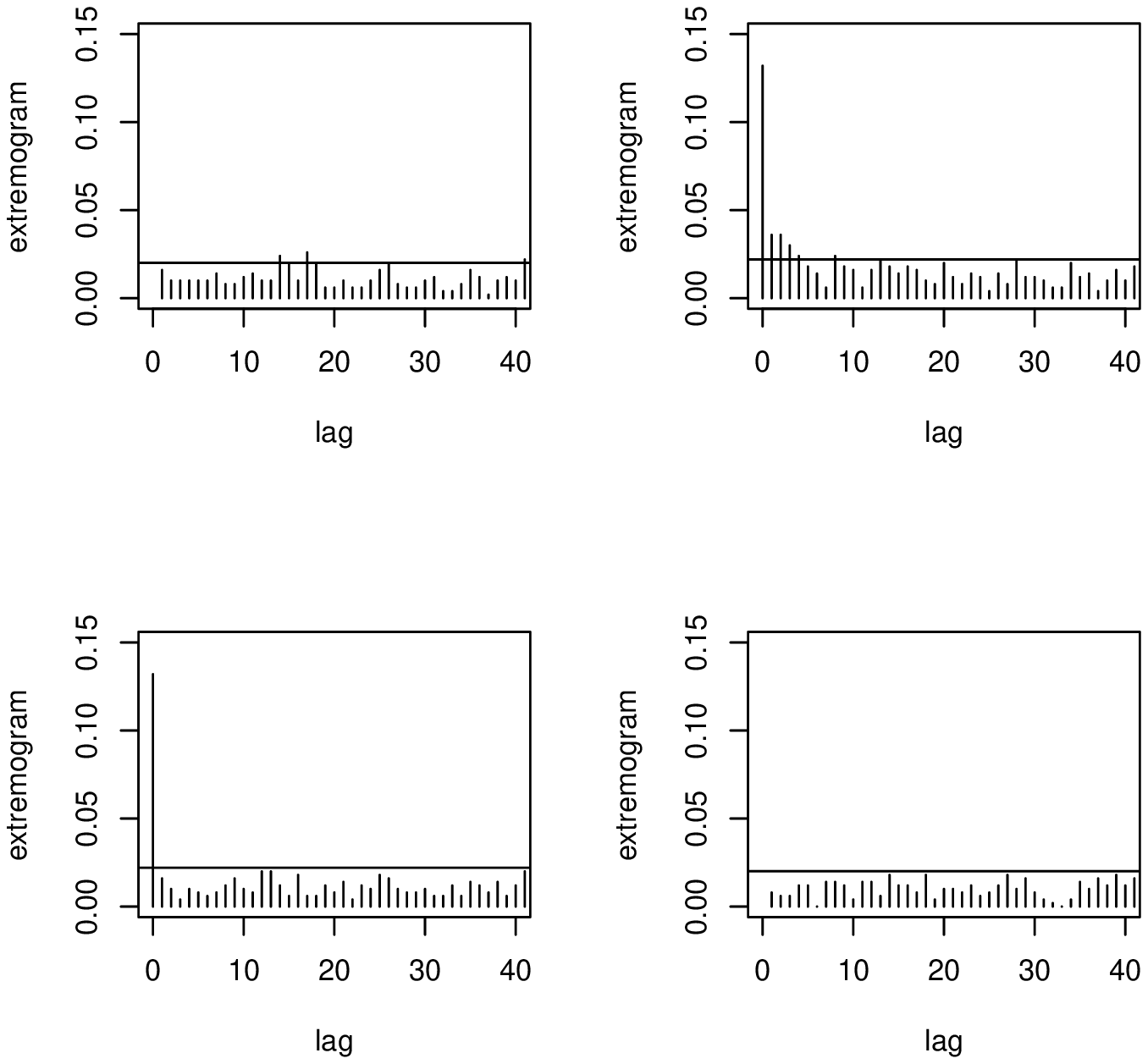}
\end{center}
\caption{(Cross-) extremograms in Example (12). Top left: bivariate \garch\ data.
Top right $2\times 2$ graphs: (Cross-) extremogram of bivariate \garch\ data.  
Bottom left $2\times 2$ graphs: (Cross-) extremogram of residuals after a bivariate QMLE fit. 
Bottom right $2\times 2$ graphs: (Cross-) extremogram of 
 residuals after applying the component-wise MLE.}
\label{fig:1-12}
\end{figure}
Figures \ref{fig:1-11} and \ref{fig:1-12} (Bottom left and right $2\times 2$ graphs) indicate that
extremal cross-dependencies are not present
in the residuals of bivariate \garch\ fits with the exception at lag 0. 
This is true for the bivariate QMLE (Figures \ref{fig:1-11} and \ref{fig:1-12}:
bottom left) but, to some extent, also for
the component-wise univariate fits (Figures \ref{fig:1-11} and  
\ref{fig:1-12}: bottom right).
However, Figure~\ref{fig:1-12} (bottom right) shows that univariate fits do not remove all
cross-dependencies from the residuals (in this case the degrees of freedom were not correctly estimated).
We experimented with distinct parameters sets close to those in
Example (12) and we also replaced
univariate MLE by univariate Gaussian QMLE. In all cases,  one cannot remove all
cross-dependencies of the residuals. Therefore bivariate \garch\ fitting
is recommended if one  believes in a bivariate \garch\ model.

\subsection{An analysis of foreign exchange rates}
We analyze a bivariate high-frequency time series, consisting of  35,135
five minute returns of  USD-DEM and USD-FRF foreign exchange rates.
Throughout this subsection we choose the $98\%$ 
component-wise sample quantiles as threshold for the sample (cross-) extremograms.
In each plot the horizontal line shows the 
$96\%$ quantile obtained from  $100$ random
permutations of the data.
\par
\begin{figure} 
\begin{center}
\includegraphics[width=0.475\textwidth]{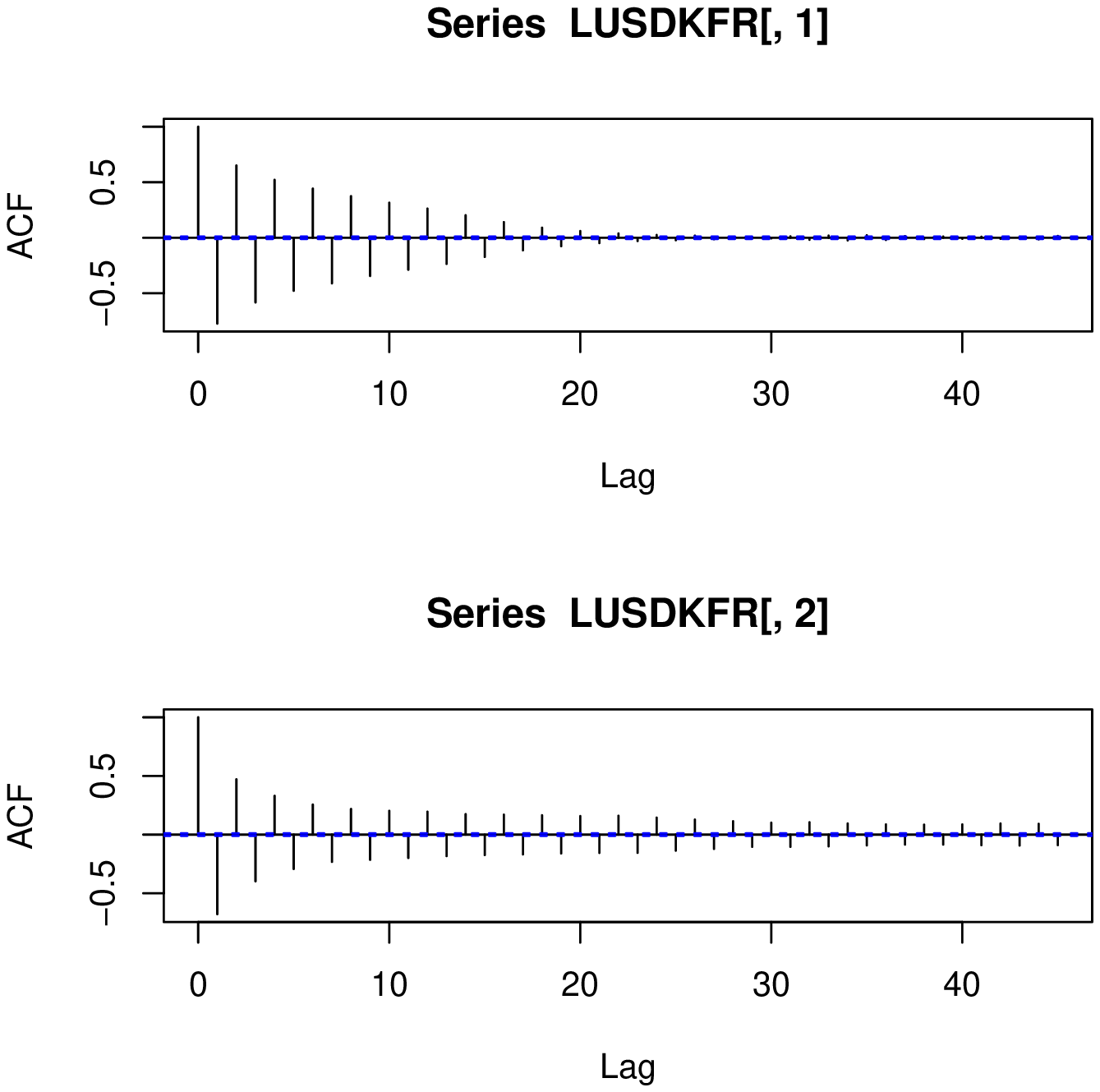}
\includegraphics[width=0.475\textwidth]{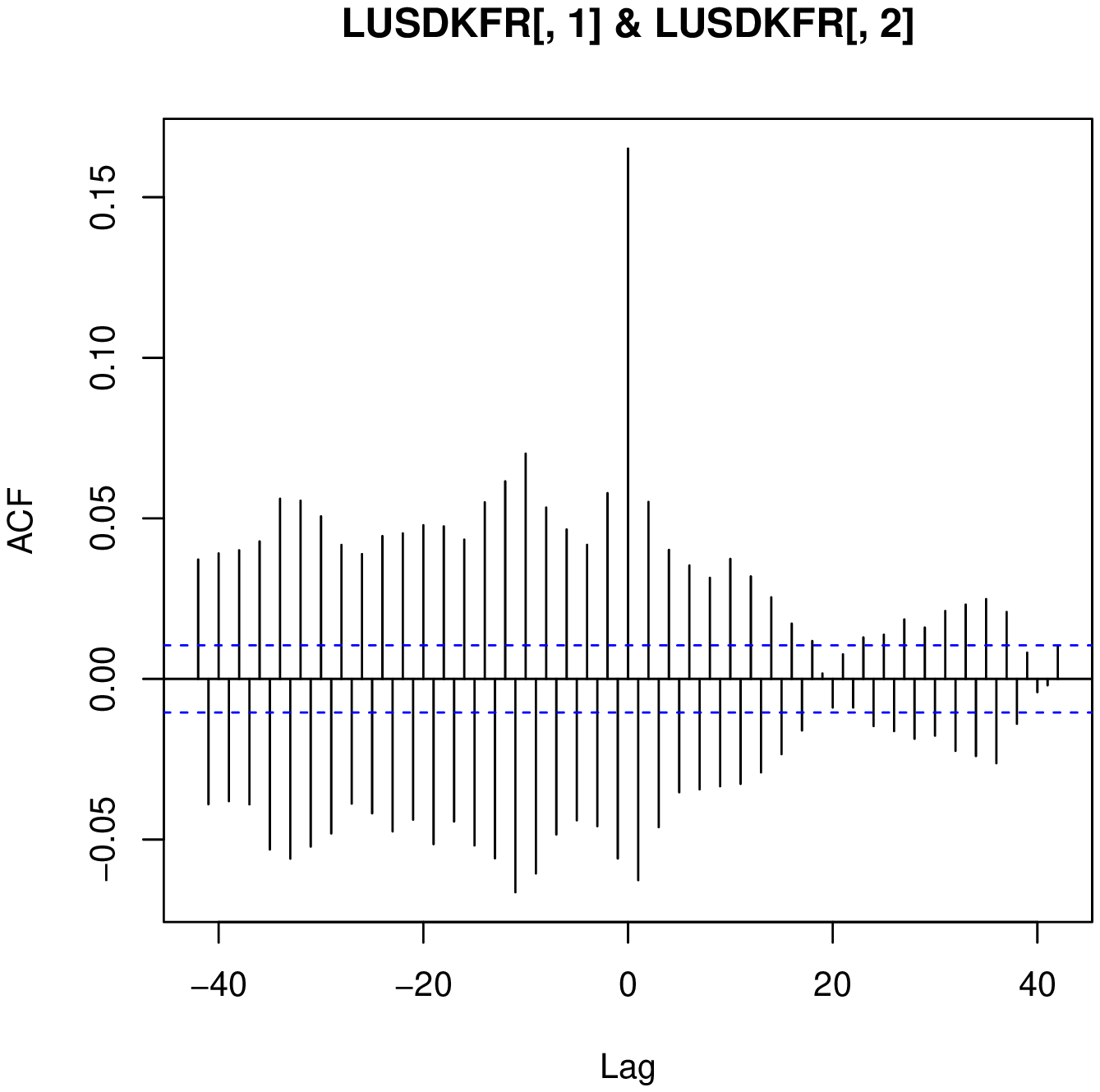} \\
\includegraphics[width=0.475\textwidth]{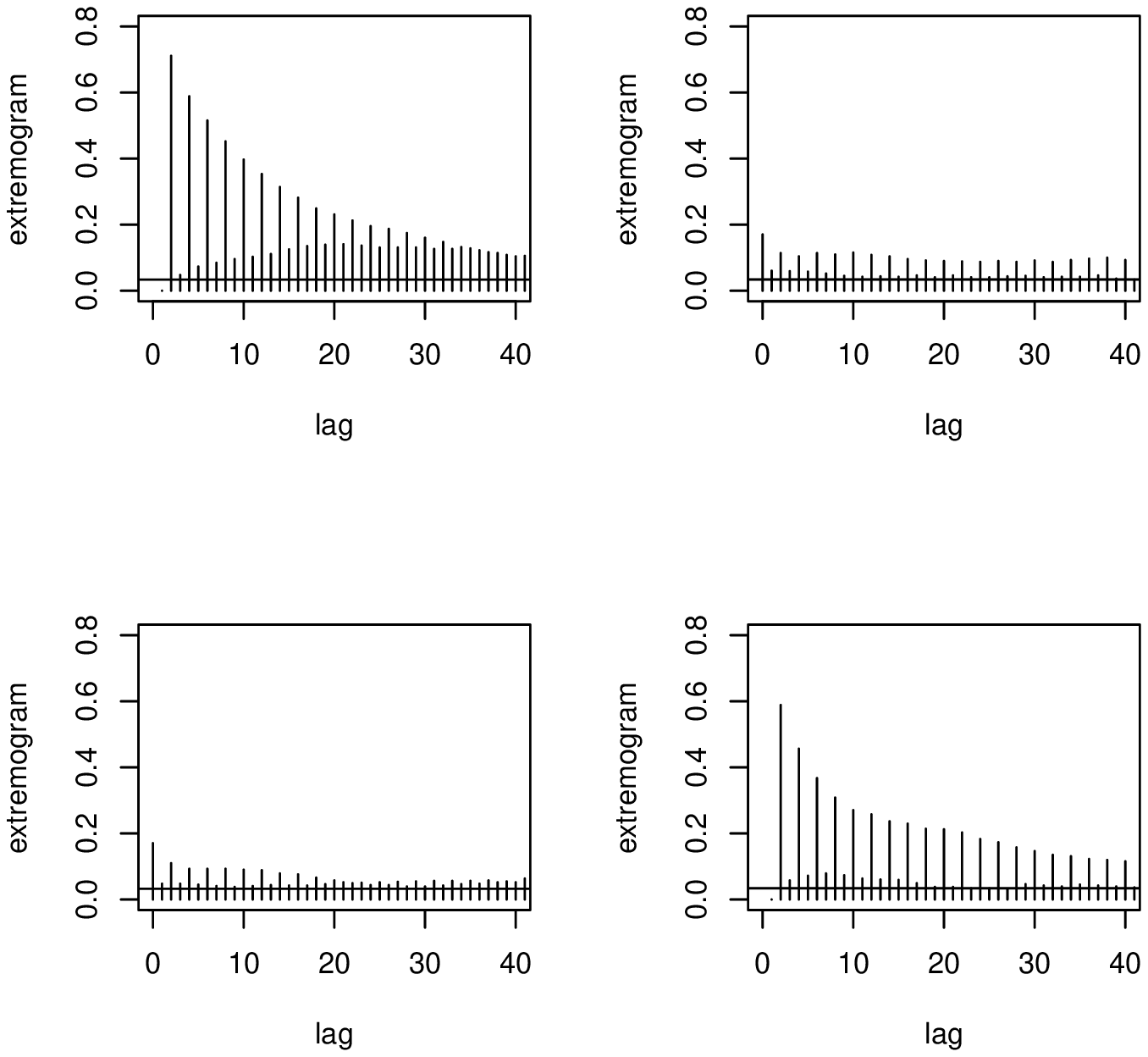}
\end{center}
\caption{Five minute log return series of USD-DEM and USD-FRF 
foreign exchange rates. 
 Top left: The sample autocorrelation \fct s for each series. Top right:
 Cross-sample autocorrelation \fct s. Bottom $2\times 2$ graphs:
 (Cross-) extremograms of the original data. (Cross-) sample
 autocorrelations oscillate strongly and show rather strong dependencies in every other lags.
 (Cross-) extremograms take over this tendency.}
\label{fig:2-1}
\end{figure}

The data exhibits rather strong 
cross-correlations and autocorrelations (see Figure \ref{fig:2-1}, 
top left: sample autocorrelations, top right: cross-correlations). 
So it is not unexpected that we also
observe dependence of the extreme values of the two series. This is apparent
in the extremograms of Figure \ref{fig:2-1} (bottom). 
After fitting a bivariate vector AR model of order $19$ to the data which is chosen by the ``Schwarz
criterion'' (``Akaike's final prediction error criterion'' proposes an order of 
$44$), we fitted a bivariate \garch\ model to the residuals, by employing
bivariate QMLE. The estimated  matrices \eqref{para-alpha-beta} are as follows:
\beao
\barr{llll}
&\left(
    \begin{array}{cc}
      .214 & .013  \\
      .110 & .223
    \end{array}
  \right),&  \left(
    \begin{array}{cc}
      .697 & .008  \\
      .280 & .663
    \end{array}
  \right),&\ .372, \earr
\eeao
which satisfy the sufficient condition for stationarity of a bivariate \garch\ model; see \eqref{eq:stationarty}. 
\begin{figure} 
\begin{center}
\includegraphics[width=0.475\textwidth]{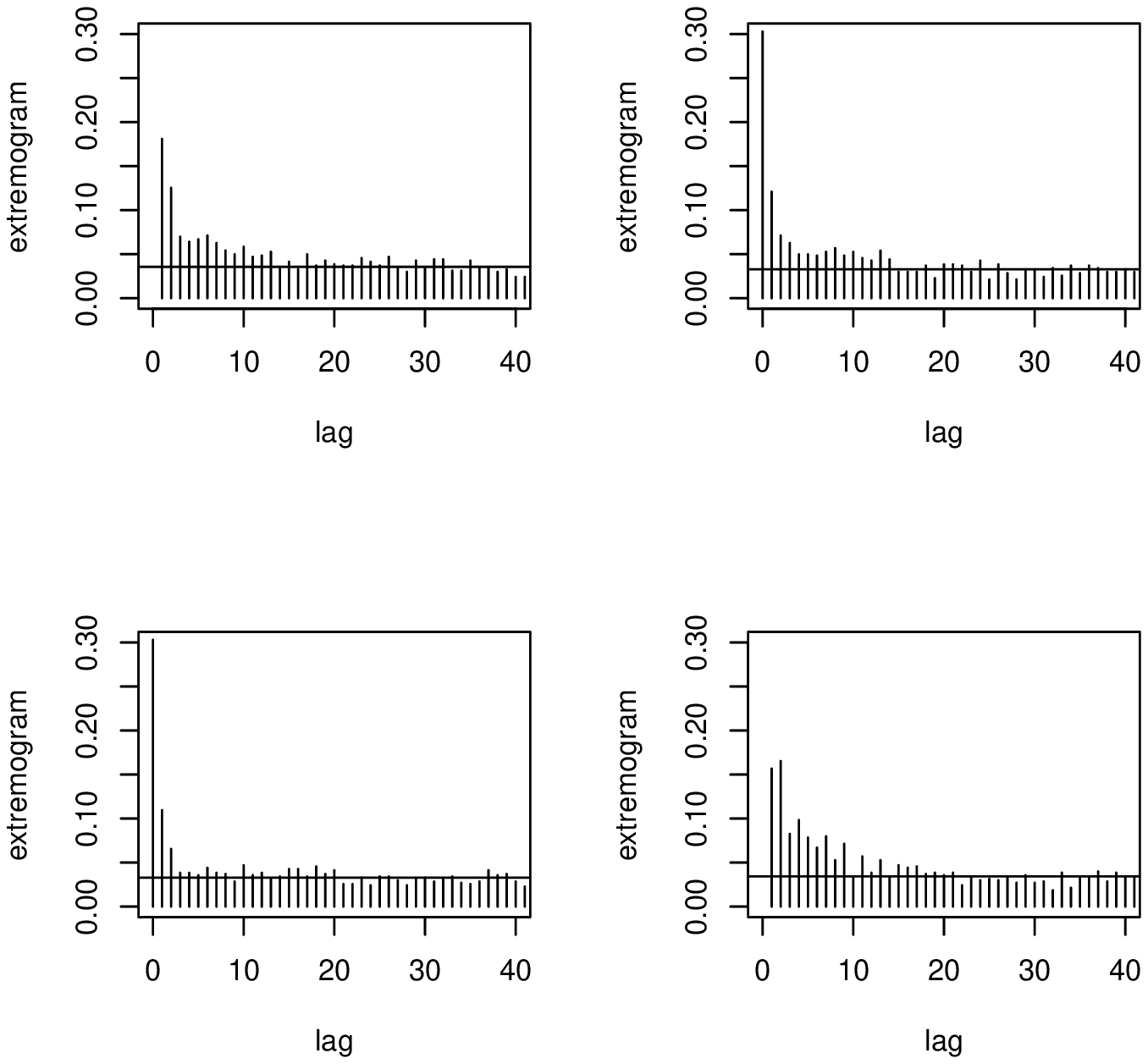}
\includegraphics[width=0.475\textwidth]{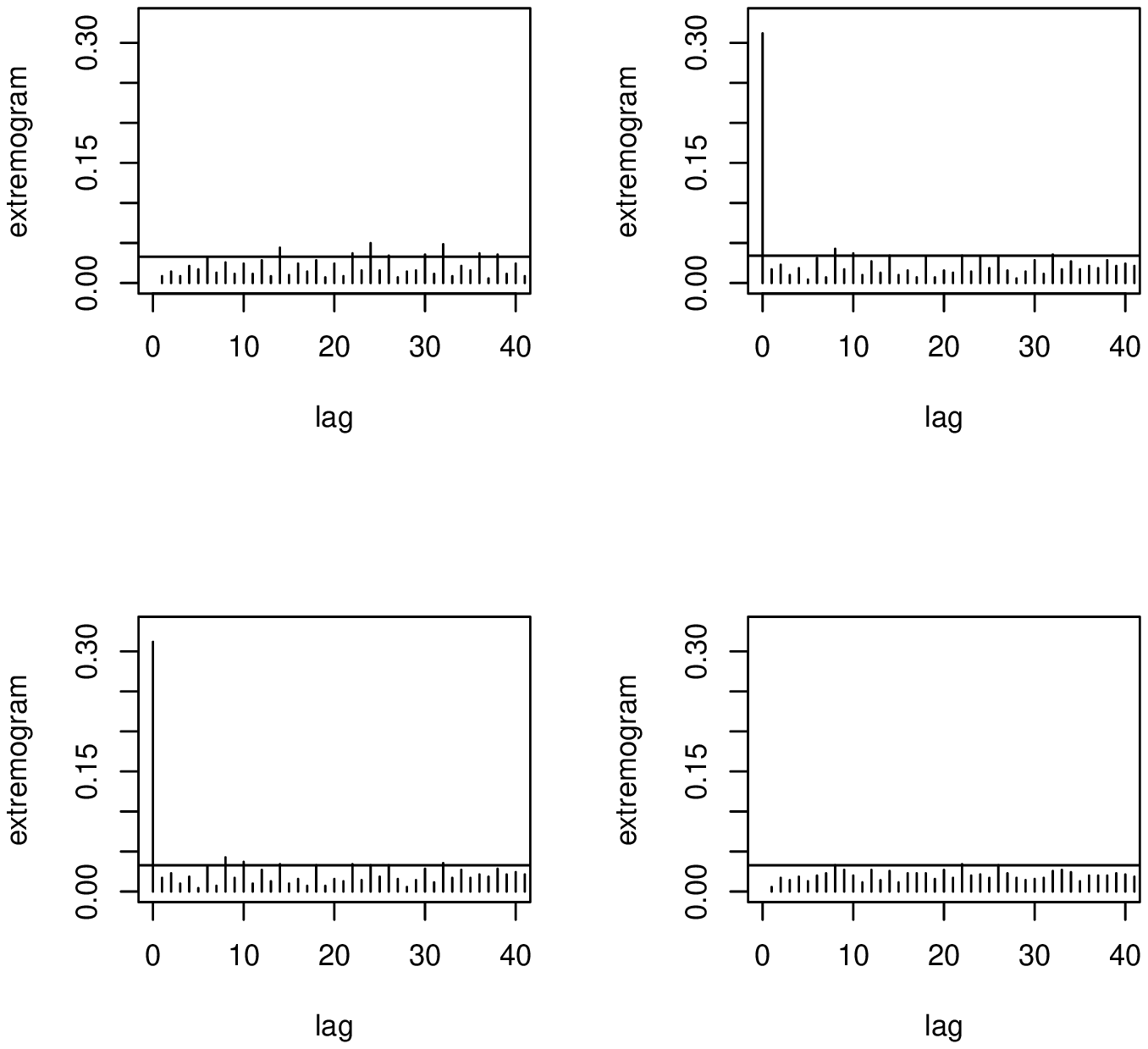}
\end{center}
\caption{(Cross-) extremograms of residuals after AR and AR-GARCH fits.
Left $2\times2$ graphs: (Cross-) extremograms of residual after AR fit. Right
 $2\times2$ graphs: Those after AR-GARCH fit. The serial (cross-) extremal dependence still remains
 in residuals of an AR fit, while it has been removed in the residuals after an AR-GARCH
 fit except for lag 0.}
\label{fig:2-3}
\end{figure}
\begin{figure} 
\begin{center}
\includegraphics[width=0.5\textwidth]{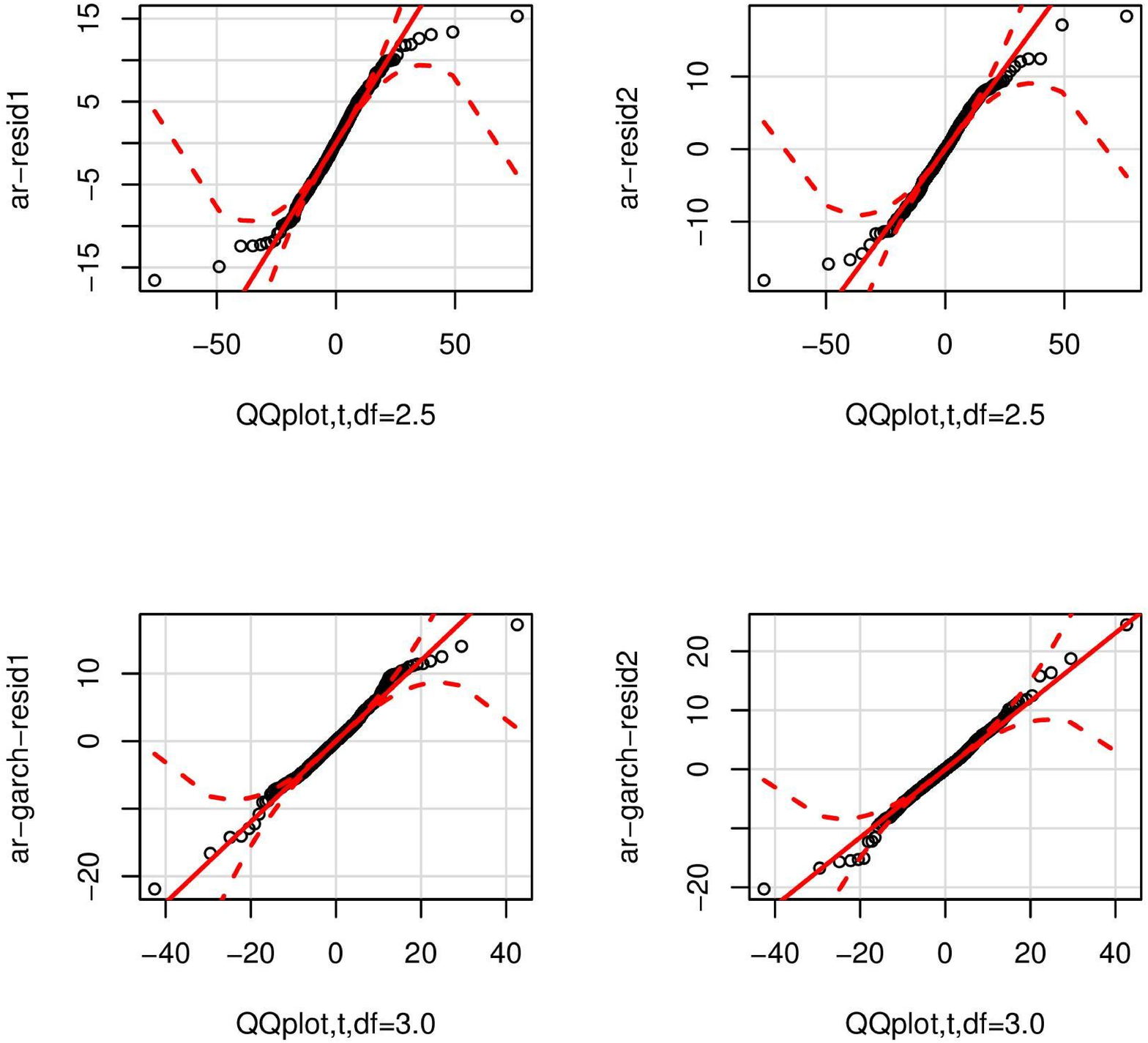}
\end{center}
\caption{QQ-plots of residuals after AR fit (top) and those after
 AR-\garch\ fit (bottom).}
\label{fig:2-4}
\end{figure}
After the AR fit,
the cross-extremograms of the residuals do  not vanish at all positive
lags; see Figure~\ref{fig:2-3} (left) and therefore a \garch\ fit for the residuals
is recommendable.
After fitting the bivariate \garch,  the  residuals exhibit
extremal cross-dependence only at time-lag $0$; see Figure~\ref{fig:2-3} (right). This means that the components of the innovations $\bfZ_t$ exhibit extremal dependence.
Judging from QQ-plots of the residuals of the vector AR model and the AR-\garch\ model
against the quantiles of $t$-\ds s with 2.5 and 3 degrees of freedom
respectively; see Figure \ref{fig:2-4}, 
these \ds s give a good fit to the \ds\ of the residuals.

\subsection{An analysis of stock returns}\label{sbsec:stock}
We consider log-return series of three stock prices from the NY Stock Exchange:
``Caterpillar Inc.'', ``FedEx Corporation'' and ``Exxon Mobil
Corporation'' (``cat'', ``fdx'' and ``xom'' for short).
In each series, the raw tick-by-tick trade data has been
processed into 5-minute grid data by taking the last realized trade price in
each interval. Prices have been restricted to exchange trading hours 9:30 a.m. to 4:00
p.m.,  Monday to Friday, so that 78 data per day have been collected in the
time period from 2009-02-18 9:30 to 2013-12-31 16:00. The
prices at 16:00 are identified with those at 9:30 of the next days.
Since we consider the log-return series, the total size of
the data is $78\times 1265+1=98671$. 
\par
\begin{figure}
\begin{center}
\includegraphics[width=0.8\textwidth]{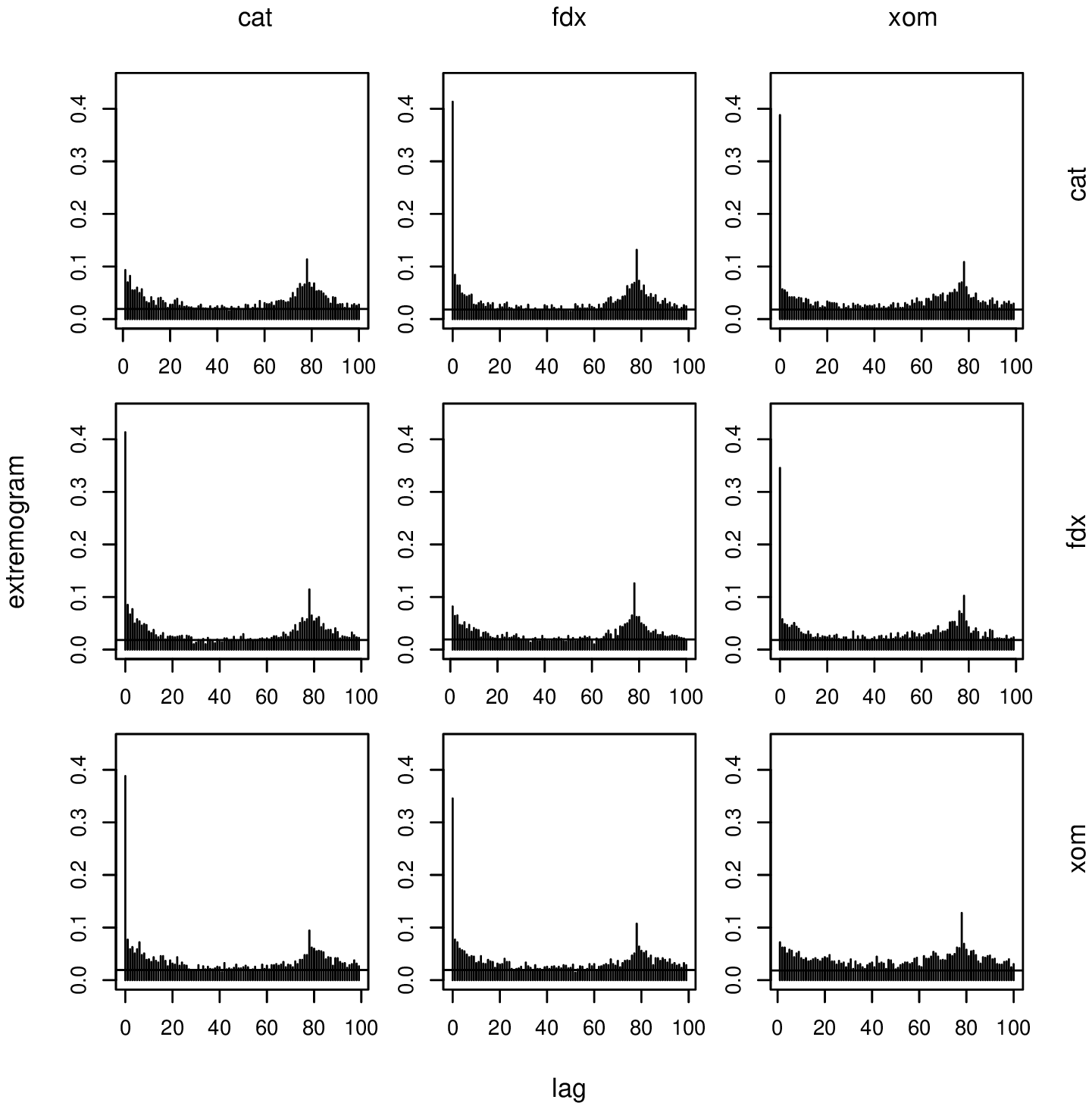}
\end{center}
\vspace{-0.5cm}
\caption{(Cross-) extremograms of log-return series of 3 stock prices (cat, fdx, xom).
 Graphs show strong serial extremal dependence in each series together
 with strong extremal dependence between $3$ series. Other than large
 spikes at lag $0$ in cross-extremograms, we observe spikes  
 at lag $78$, which show seasonal fluctuation 
 in a day. Moreover, extremal data around the begging 9:30 a.m. and the
 end 4:00 p.m. may exhibit strong dependence. \label{fig:3-1}}
\end{figure}
\begin{figure}
\begin{center}
\includegraphics[width=0.9\textwidth]{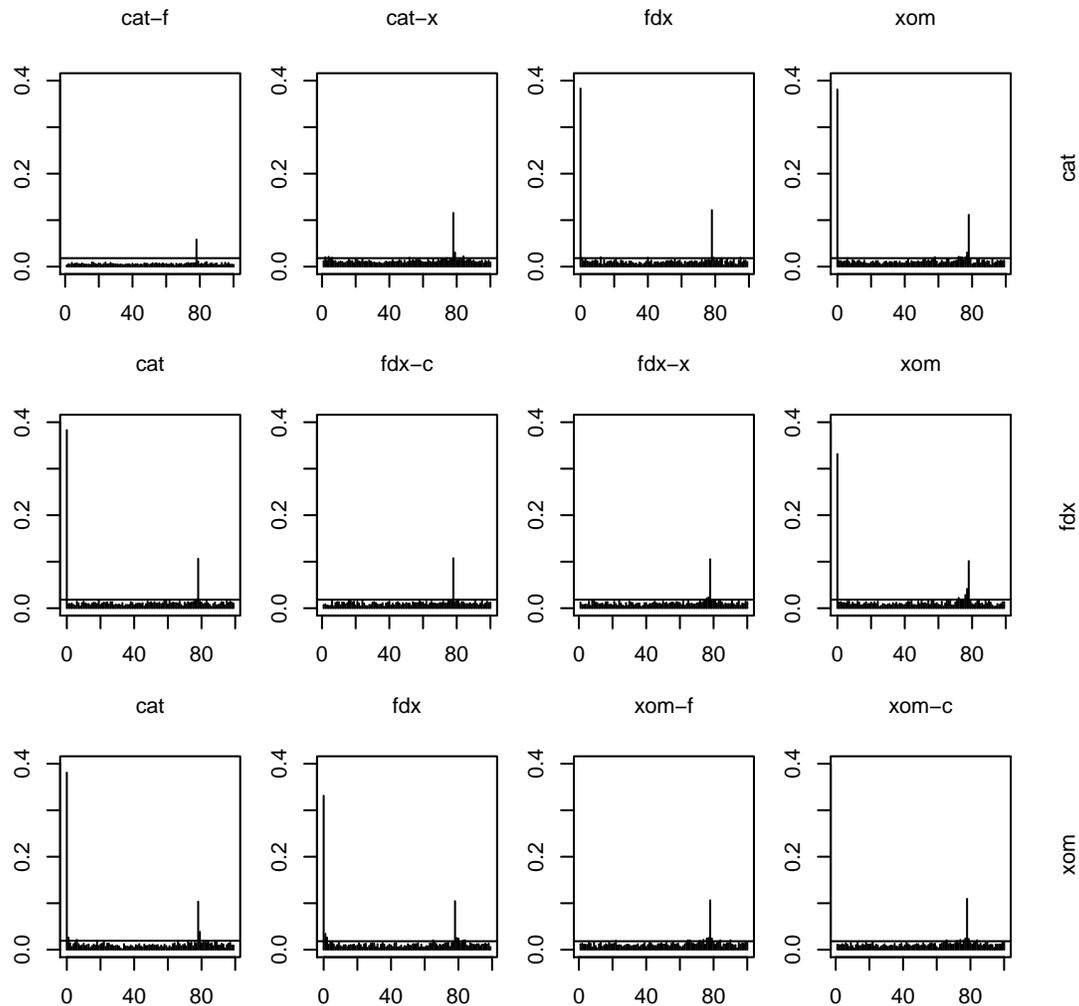}
\vspace{-1cm}
\caption{(Cross-) extremograms for residuals of bivariate
 \garch\ fits to combinations (cat, fdx), (fdx, xom) and (xom, cat), so that 
we have two extremograms in each row, i.e. cat-f and cat-x are those for
 residuals of 'cat' components respectively from bivariate QMLE of (cat,
 fdx) and (xom, cat). Other elements are cross-extremograms
 for residuals of (cat, fdx, xom) against (cat, fdx, xom). Although,
 residuals show less extremal dependence except for large spikes at
 $0$, we could not remove the seasonal component at lag $78$. \label{fig:3-2}}
\ece
\end{figure}

The sample (cross-)~extremograms of the log-returns of the stock prices are shown in Figure
\ref{fig:3-1}, where we choose the empirical $0.99$ quantiles of the returns
as the threshold. We take the 
$97\%$ quantile obtained from  $100$ random
permutations of the data and show the values as the horizontal line
in each plot. 
Although 
we observe typical \garch\ (cross-) extremograms close to lag 
0, there is a clear seasonal component in these plots, appearing as spikes 
at lag 78, 
corresponding to the beginning and end of the days. 
A \garch\ model (bivariate or component-wise univariate) cannot explain the seasonal extremal components
in the data. However, the (cross-) extremograms of the residuals after a bivariate \garch\ fit show that most of the
serial dependence has been removed from the data, although the seasonal component is also present in the residuals;  
see Figure~\ref{fig:3-2}. 


We fit a bivariate \garch\ model to each pair of stock prices,
i.e., (cat,fdx), (fdx,xom) and (cat,xom). The estimated 
values of the bivariate QMLE \eqref{para-alpha-beta} for (cat,fdx), (fdx,xom), (cat,xom), respectively, are 
\begin{align*}
& \left(
    \begin{array}{cc}
      .215 & .210  \\
      .029 & .287
    \end{array}
  \right),  \left(
    \begin{array}{cc}
      .666 & .144  \\
      .002 & .668
    \end{array}
  \right),\ .55,  \quad 
\left(
    \begin{array}{cc}
      .178 & .000  \\
      .006 & .250
    \end{array}
  \right), \left(
    \begin{array}{cc}
      .712 & .115  \\
      .007 & .666
    \end{array}
  \right),\ .484,  \\
& \hspace{3cm}\left(
    \begin{array}{cc}
      .094 & .153  \\
      .009 & .278
    \end{array}
  \right),  \left(
    \begin{array}{cc}
      .789 & .094  \\
      .007 & .650
    \end{array}
  \right),\ .567.
\end{align*}
The estimators of the combination (cat, fdx) are unstable and take
the boundary value\footnote{
In this case, the univariate \garch\ fit does not converge. Therefore we
examine several initial values for ``ccgarch'' on a grid of size $0.1$
and choose an ``optimal'' value based on their likelihoods. We also tried
several optimization methods included in ``ccgarch''.  Then we calculated
the  eigenvalues of \eqref{eq:spectral} from the estimates, including the ``optimal'' ones. 
However, the largest eigenvalues
 are very close to one  in all cases. Since ``ccgarch'' finds the
 optimal value under the sufficient condition  \eqref{eq:stationarty},  the real optima
 would certainly violate \eqref{eq:stationarty}. }
of the sufficient condition \eqref{eq:stationarty} while the  
 estimates for (fdx, xom) and (cat, xom) satisfy \eqref{eq:stationarty}. The
obvious seasonal component  of the data 
(corresponding
 to the end of a trading day at lag 78)
probably violates the stationary
 condition.  
Nevertheless, the standardized residuals appear ``de-volatilized'' in all
(cross-) extremograms modulo the fact that the seasonal component is always present.
 
\begin{figure}
\begin{center}
\includegraphics[width=0.8\textwidth]{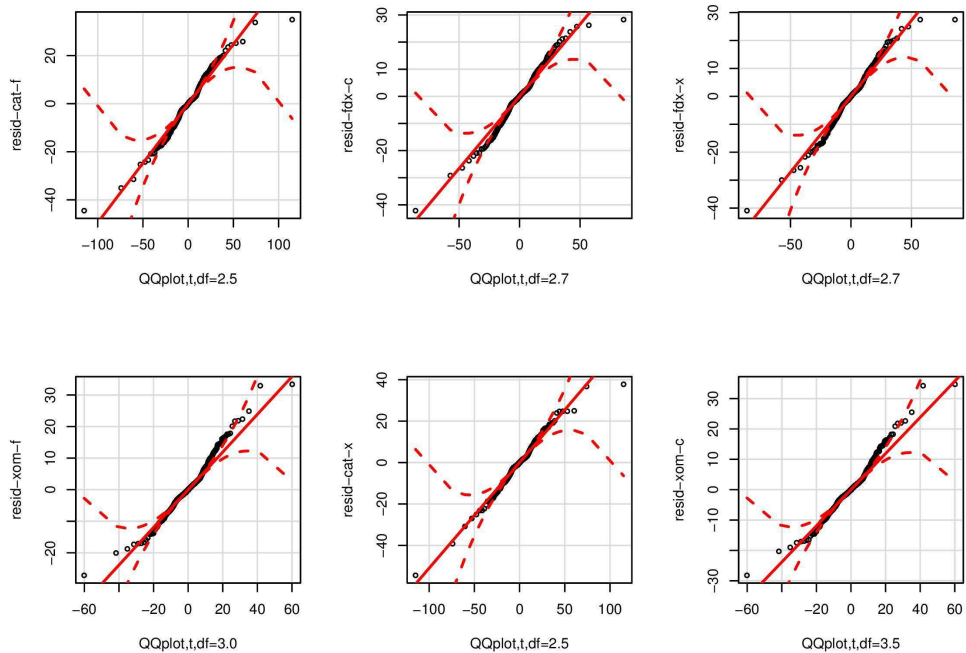}
\caption{QQ-plots of the residuals after the bivariate \garch\ fits against 
the quantiles of $t$-\ds s. Top left and middle: QQ-plots for residuals of
 'cat' and 'fdx' components respectively, after the bivariate
 \garch\ fits to the pair (cat, fdx). Top right and bottom left: Those for
 residuals of 'fdx' and 'xom'. Bottom left and right: Those for
 residuals of 'cat' and 'xom'.  
} \label{fig:3-4}
\end{center}
\end{figure}

We take QQ-plots of the residuals after bivariate \garch\ fits against 
the quantiles of $t$-\ds\ with $2.5 \sim 3.5$ degrees of freedom
respectively; see Figure \ref{fig:3-4}. From these plots, $t$-distributions show
a good fit to the \ds\ of the residuals, which also assures that innovations of
 \garch\ model satisfy the regular variation assumption. Only residuals
 of 'xom' components (bottom left and bottom right in Figure
 \ref{fig:3-4}) seem to fit $t$-\ds\ with different degrees, depending on
 the pair of stock prices in bivariate estimations. 

\begin{figure}
\begin{center}
\includegraphics[width=0.3\textwidth]{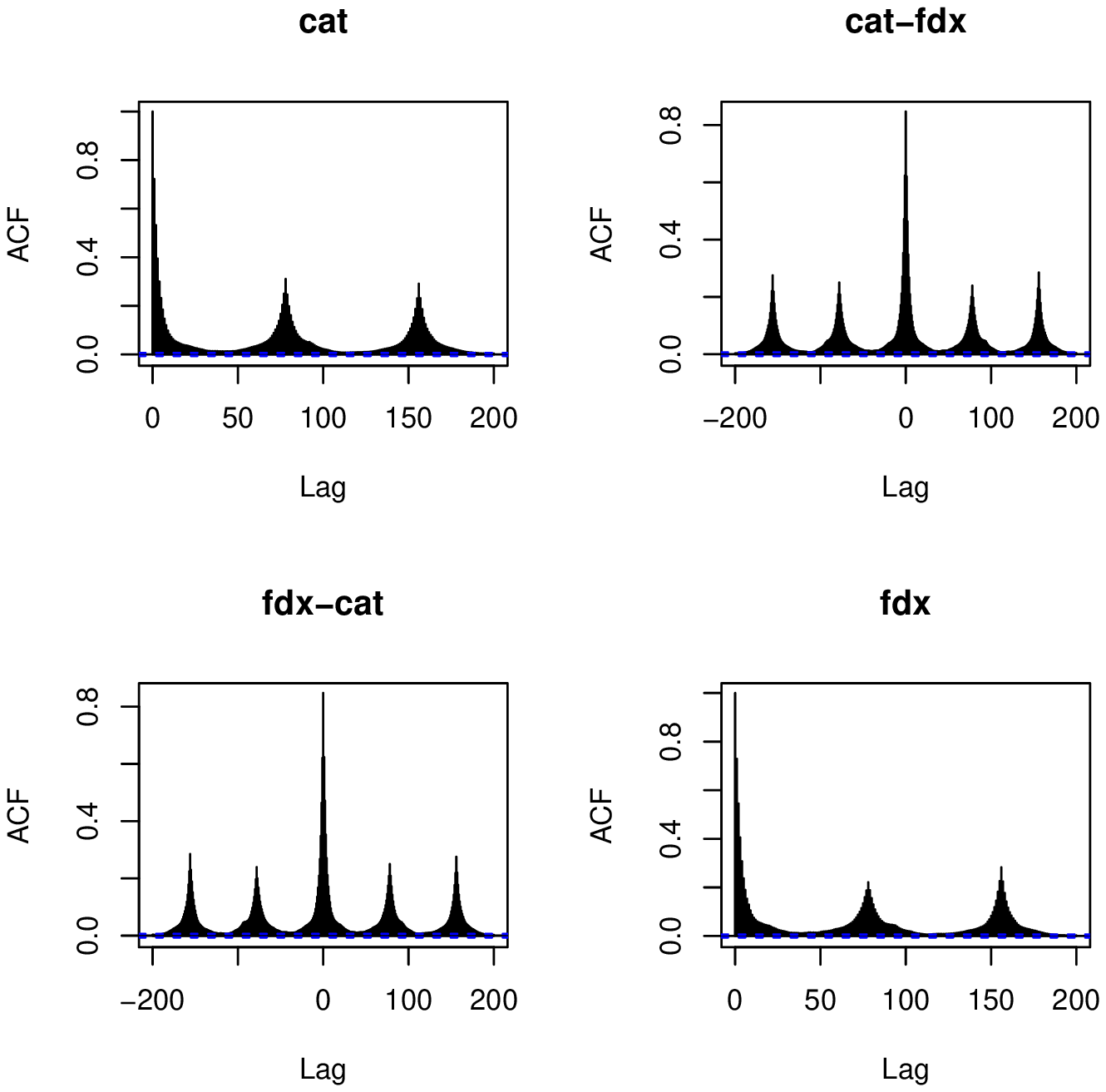}
\includegraphics[width=0.3\textwidth]{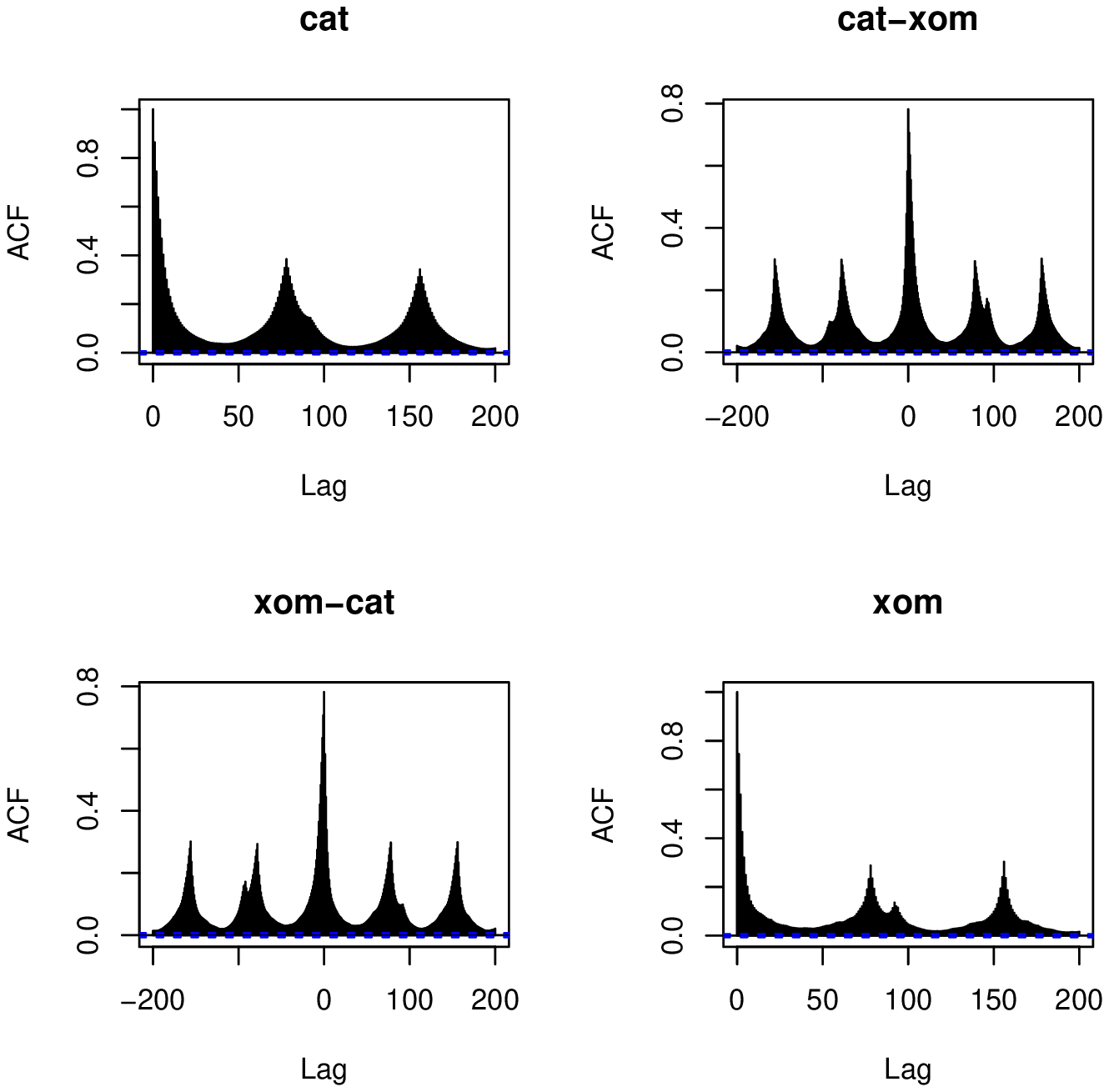}
\includegraphics[width=0.3\textwidth]{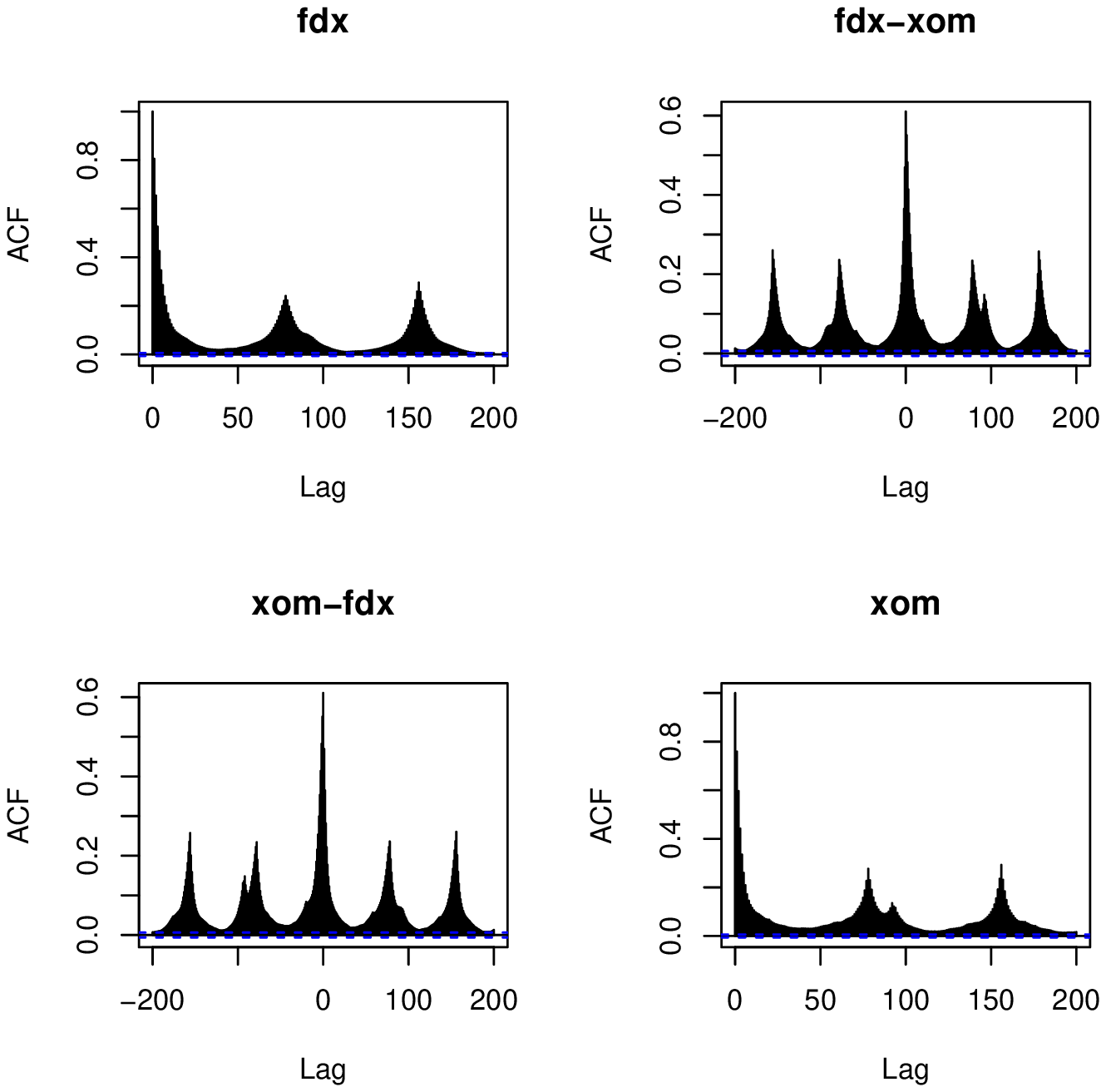}
\caption{(Cross-) autocorrelation functions for estimated volatility processes by bivariate
 \garch\ fits. Left $2\times 2$ graphs: (Cross-) autocorrelation functions for estimated volatilities after bivariate
 \garch\ fits to the pair (cat, fdx). Middle $2\times 2$ graphs: Those for
 the pair (fdx, xom). Right $2\times 2$ graphs: Those for the pair (xom, cat). 
 In all (cross-) autocorrelation functions, one could observe clear 
 seasonal effects by a day cycle ( we have 78
 of 5-minute intervals per day from 9:30 to 16:00.). \label{fig:3-3}}
\end{center}
\end{figure}

As predicted, estimated volatilities are
affected by periodicity of the day, which is also approved in (cross-)
 autocorrelation functions of estimated volatility processes in Figure
 \ref{fig:3-3}. Since the seasonal patterns in these plots are quite
 clear and intuitive, it is desirable to remove the seasonal component
 for further analysis. 

\begin{figure}
\begin{center}
\includegraphics[width=0.6\textwidth]{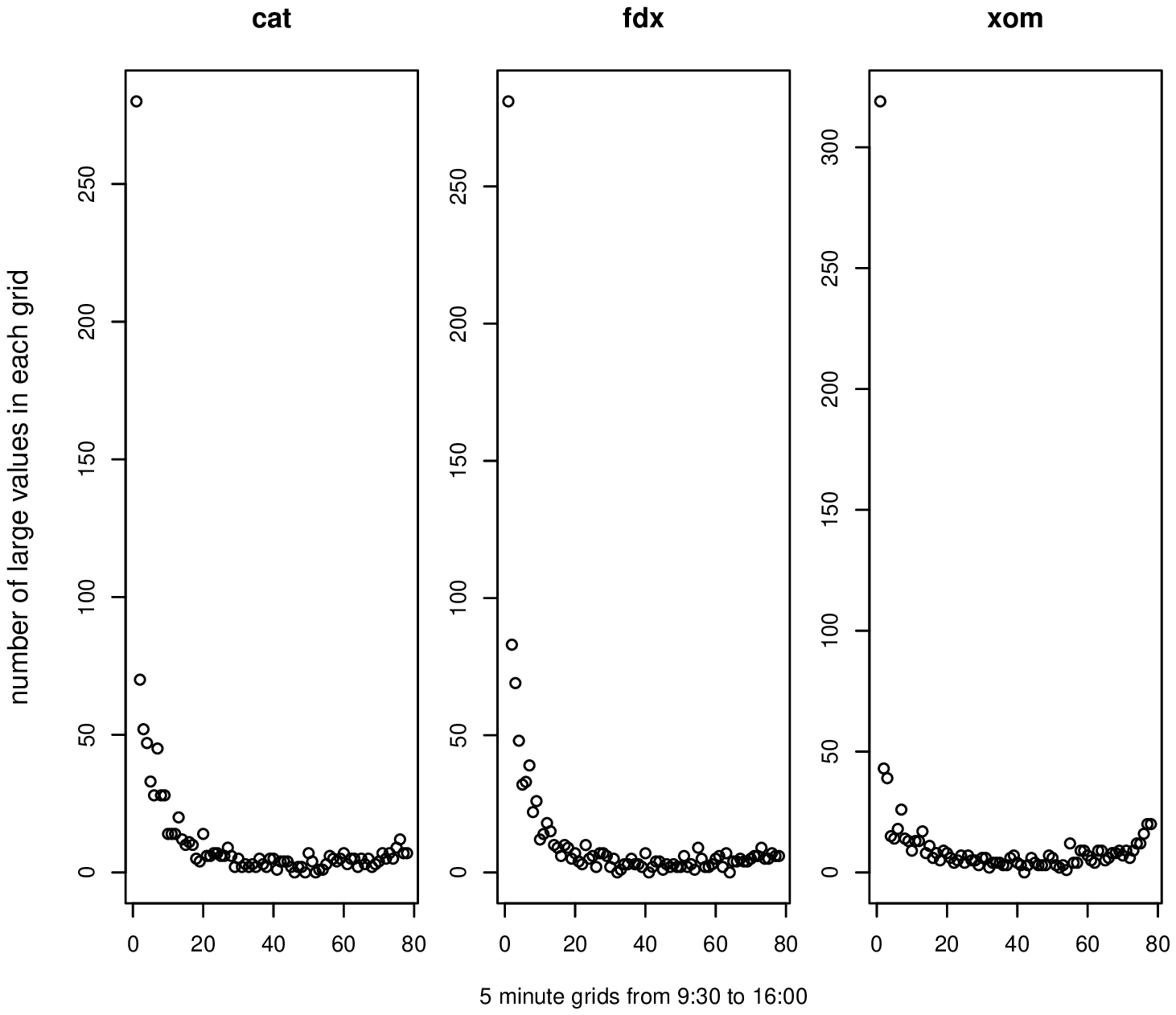}
\caption{Frequency of large values of log-return series in each 5-minute
 grid in a day. There are 78 grids in a day and we chose values larger
 then the upper $0.99$ empirical quantiles of $78 \times 1265=98670$ data. In all cases, large values
 concentrate around the beginning and end of exchange trading hours. } \label{fig:3-5}
\end{center}
\end{figure}

We further investigated the seasonal effects by plotting the number of 
large values in each 5-minute grid during the exchange trading hours; see Figure
 \ref{fig:3-5}. The large values tend to appear around the begging and end of the
 trading hours. This is typically observed in the log-return series of stock 
 price data. \\

\noindent {\bf Acknowledgment.}
We would like to thank Martin Anders J\"onsson for arranging us the
 stock price data. A major part of this work was done when Muneya Matsui
 visited the Department of Mathematics at the University of
 Copenhagen. He is grateful to its hospitality.

\end{document}